\documentclass[]{amsart}
\usepackage{amsmath}
\usepackage{amssymb,amsfonts}
\usepackage{xypic}
\xyoption{all}

\newtheorem{theorem}{Theorem}[section]
\newtheorem{definition}[theorem]{Definition}
\newtheorem{proposition}[theorem]{Proposition}

\newtheorem{lemma}[theorem]{Lemma}
\newtheorem{remark}[theorem]{Remark}

\def\m{\ensuremath{\mathcal{M}}}

\def\a{\ensuremath{\mathcal{A}}}
\def\b{\ensuremath{\mathcal{B}}}
\def\c{\ensuremath{\mathcal{C}}}
\def\d{\ensuremath{\mathcal{D}}}
\def\cat{\ensuremath{\mathbf{Cat}}}
\def\ner{\ensuremath{\mathrm{N}}}
\def\gner{\ensuremath{\Delta}}

\def\y{\ensuremath{\mathbf{y}}}
\def\set{\ensuremath{\mathbf{Set}}}
\def\cgner{\ensuremath{\underline{\Delta}}}
\def\cner{\ensuremath{\underline{\mathrm{N}}}}
\def\diag{\ensuremath{\mathrm{diag}}}
\newcommand{\sset}{\ensuremath{\mathbf{Simpl.Set}}}
\newcommand{\bicat}{\ensuremath{\mathbf{Bicat}}}
\newcommand{\bsset}{\ensuremath{\mathbf{Bisimpl.Set}}}
\newcommand{\class}{\ensuremath{\mathrm{B}}}

\newcommand{\lfunc}{\ensuremath{\mathrm{LaxFunc}}}

\newcommand{\clfunc}{\ensuremath{\underline{\mathrm{LaxFunc}}}}
\newcommand{\colfunc}{\ensuremath{\underline{\mathrm{OpLaxFunc}}}}
\newcommand{\cnlfunc}{\ensuremath{\underline{\mathrm{NorLaxFunc}}}}
\newcommand{\conlfunc}{\ensuremath{\underline{\mathrm{NorOpLaxFunc}}}}
\newcommand{\nhom}{\ensuremath{\underline{\mathrm{NorHom}}}}
\newcommand{\hoco}{\ensuremath{\mathrm{hocolim}}}
\newcommand{\func}{\ensuremath{\mathrm{Func}}}

\def\xlar(#1){\xymatrix@1{ & \ar[l]_{#1}}}
\def\xar(#1){\xymatrix@1{ \ar[r]^{#1}& }}
\begin{document}
\title[{\em Classifying Spaces for Bicategories}]
{{Nerves and Classifying Spaces for Bicategories}}

\author{P. Carrasco}
\author{A.M. Cegarra}
\author{A. R. Garz\'{o}n}
\thanks{ Partially supported by DGI of Spain and FEDER (Project: MTM2007-65431); Consejer\'{i}a de Innovacion
de J. de Andaluc\'{i}a (P06-FQM-1889); MEC de Espa\~{n}a, `Ingenio Mathematica(i-Math)' No.
CSD2006-00032 (consolider-Ingenio 2010).}
\address{
Departamento de \'Algebra \newline Facultad de Ciencias \newline Universidad de Granada
\newline 18071 Granada, Spain \newline mcarrasc@ugr.es\newline
acegarra@ugr.es\newline agarzon@ugr.es}

\begin{abstract}
This paper explores the relationship amongst the various simplicial and pseudo-simplicial objects
characteristically associated to any bicategory $\c$. It proves the fact that the geometric
realizations of all of these possible candidate `nerves of $\c$' are homotopy equivalent. Any one
of these realizations could therefore be taken as {\em the} classifying space $\class\c$ of the
bicategory. Its other major result  proves a direct extension of Thomason's `Homotopy Colimit
Theorem'  to bicategories: When the homotopy colimit construction is carried out on a diagram of
spaces obtained by applying the classifying space functor to a diagram of bicategories, the
resulting space has the homotopy type of a certain bicategory, called the `Grothendieck
construction on the diagram'. Our results provide coherence for all reasonable extensions to
bicategories of Quillen's definition of the `classifying space' of a category as  the geometric
realization of the category's Grothendieck nerve, and they are applied to monoidal (tensor)
categories through the elemental `delooping' construction.

\end{abstract}


\maketitle


\section{Introduction and summary}
Higher-dimensional categories provide a suitable setting for the treatment of an extensive list of
subjects with recognized mathematical interest. The construction of nerves and classifying spaces
of higher categorical structures, and bicategories in particular, discovers ways to transport
categorical coherence to homotopical coherence  and it has shown its relevance as a tool in
algebraic topology, algebraic geometry, algebraic $K$-theory, string theory,  conformal fields
theory, and in the study of geometric structures on low-dimensional manifolds.

This paper explores the relationship amongst the various simplicial and pseudo-simplicial objects
that have been (or might reasonably be) functorially and characteristically associated to any
bicategory $\c$. It outlines or proves in detail the far from obvious fact that the geometric
realizations of all of these possible candidate `nerves of $\c$' are homotopy equivalent. Any one
of these realizations could therefore be taken as {\em the} classifying space $\class\c$ of the
bicategory. Since this result quite reasonably extends, to bicategories,  Quillen's definition
\cite{quillen} of the `classifying space' of a category such as $\class\ner\c$, the geometric
realization of the category's Grothendieck nerve, and given that monoidal (tensor) categories may
be identified with bicategories having a single object, the paper may quite possibly be of special
interest to $K$-theorists as well as to researchers interested in homotopy theory of higher
categorical structures. Moreover, the other non-elementary result of the paper gives and
non-trivially proves a direct extension of Thomason's `Homotopy Colimit Theorem' \cite{thomason} to
bicategories: When the homotopy colimit construction is carried out on a diagram of spaces obtained
by applying the classifying space functor to a diagram of bicategories, the resulting space has the
homotopy type of a certain bicategory called the `Grothendieck construction on the diagram'.
Notice that, even in the case where the diagram is of monoidal categories and monoidal functors,
the Grothendieck construction on it turns out to be a genuine bicategory. Hence, the reader
interested in the study of classifying spaces of monoidal categories can find in the above fact a
good reason to also be interested in the study of classifying spaces of bicategories.

There is a miscellaneous collection of ten different `nerves'  for a bicategory $\c$ discussed in
the paper, each with a particular functorial property.  All these nerves occur in a commutative
diagram
$$
\xymatrix{&&&\ \  \cner\c \ar@{~>}@<1.5pt>[d]\ar@{~>}@<1.5pt>[dl]\ \ \ar@{~>}@<-1.5pt>[dr]&&\\
\mathrm{(A)} &\cgner\c\ & \ar@{->}@<1.5pt>[l]\ \cgner\hspace{-4pt}^{^\mathrm{u}}\c\ &
\ar@{->}@<1.5pt>[l]\ \underline{\mathrm{S}}\,\c \ \ar@{->}@<-1.5pt>[r] &\
\overline{\nabla}\hspace{-2.5pt}_{^{_\mathrm{u}}}\c\  \ar@{->}@<-1.5pt>[r]& \ \overline{\nabla}\c
\\ & \ar@{->}@<1.5pt>[u]\gner\c\
&\ar@{->}@<1.5pt>[l]\ \gner\hspace{-4pt}^{^\mathrm{u}}\c \ \ar@{->}@<1.5pt>[u] &&\
\nabla\hspace{-2.5pt}_{^{_\mathrm{u}}}\c \ \ar@{->}@<1.5pt>[u] \ar@{->}@<-1.5pt>[r] & \ \nabla\c, \
\ar@{->}@<1.5pt>[u]}
$$
in which the arrows written as $\to$ denote simplicial maps (those in the bottom row) or simplicial
functors (all the others),  whereas those written as $\rightsquigarrow$ are pseudo-simplicial
functors.

In the diagram, $\cner\c:\Delta^{\!{^\mathrm{op}}} \rightsquigarrow \cat$ is the pseudo-simplicial
category whose category of $p$-simplices is
$$
\cner\c_p = \bigsqcup_{(x_0,\ldots,x_p)\in \mbox{\scriptsize Ob}\c^{p+1}}\hspace{-0.3cm}
\c(x_1,x_0)\times\c(x_2,x_1)\times\cdots\times\c(x_p,x_{p-1}),
$$
where $\c(x,y)$ denotes the hom-category of the bicategory $\c$ at a pair of objects $(x,y)$. The
face and degeneracy functors are defined in the standard way by using the horizontal composition
and the identity morphisms of the bicategory,  and the natural isomorphisms $d_id_j\cong
d_{j-1}d_i$, etc., being given from the associativity and unit constraints of the bicategory.
Therefore, when a category $\c$ is considered as a discrete bicategory, that is, where the
deformations (2-cells) are all identities, then $\cner\c=\ner \c$, the usual Grothendieck's nerve
of the category. This $\cat$-valued `pseudo-simplicial nerve' of a bicategory is explained in some
detail in the third section of the paper, and its realization is taken as the `classifying space'
$\class \c$ of the bicategory. Since the horizontal composition involved is in general neither
strictly associative nor unitary, $\cner\c$ is not a simplicial category. Consequently, we are,
unfortunately, forced to deal with defining the geometric realization of what is not simplicial but
only `simplicial up to isomorphisms'. Indeed, as we review in the preliminary Section 2, this has
been done by Segal, Street and Thomason using some Grothendieck's methods, but the process is quite
indirect and the $CW$-complex $\class\c$ thus obtained has little apparent intuitive connection
with the cells of the original bicategory $\c$.

However, the rest of the nerves in the diagram above do not have these simplicial defects caused by
the lack of associativity  or unitary properties of the horizontal composition in the bicategory.
To refer briefly to some of them here, say that, for example,  $\gner\hspace{-4pt}^{^\mathrm{u}}\c:
\Delta^{\!{^\mathrm{op}}} \to \set$ and $\nabla\hspace{-2.5pt}_{^{_\mathrm{u}}}\c:
\Delta^{\!{^\mathrm{op}}} \to \set$ are the nerves of the bicategory introduced by Street in
\cite{street} and Duskin in \cite{duskin}. These are genuine single simplicial sets, termed here
the `unitary geometric nerves' of the bicategory. As observed by Duskin (see also \cite{gurski2},
for an interesting new approach), both $\gner\hspace{-4pt}^{^\mathrm{u}}\c$ and
$\nabla\hspace{-2.5pt}_{^{_\mathrm{u}}}\c$ completely encode all the structure of the bicategory
and, furthermore, there is a pleasing geometrical description of their simplices: a $p$-simplex of
$\Delta^{\hspace{-2pt}^\mathrm{u}}\c$ (respect. of $\nabla\hspace{-2.5pt}_{^{_\mathrm{u}}}\c$) is
geometrically represented by a diagram in $\c$ with the shape of the 2-skeleton of an oriented
affine standard $p$-simplex, whose faces are triangles
$$
\xymatrix{\ar@{}[drr]|(.6)*+{\Downarrow \widehat{F}_{_{i,j,k}}}
                & Fj \ar[dl]_{\textstyle F_{i,j}}             \\
Fi  & &     Fk\,  \ar[ul]_{\textstyle F_{j,k}} \ar[ll]^{\textstyle F_{i,k}} }\hspace{1cm} \xymatrix{ \ar@{-}@/_0.5pc/[d]&\ar@{}[drr]|(.6)*+{\Uparrow \widehat{F}_{_{i,j,k}}} & Fj \ar[dl]_{\textstyle \text{respect. }\hspace{0.2cm}\textstyle F_{i,j}}         \\
 & Fi  & &     Fk   \ar[ul]_{\textstyle F_{j,k}} \ar[ll]^{\textstyle F_{i,k}} }  \xymatrix@R=13pt{\ar@{-}@/^0.5pc/[dd]&\\&\\&}
$$
with objects (0-cells) $Fi$  placed on the vertices, $0\leq i\leq p$, morphisms (1-cells)
$F_{i,j}:Fj\rightarrow Fi$ on the edges, $0\leq i<j\leq p$, and $\widehat{F}_{i,j,k}:F_{i,j}\circ
F_{j,k}\Rightarrow F_{i,k}$ (respect. $\widehat{F}_{i,j,k}:F_{i,k}\Rightarrow F_{i,j}\circ F_{j,k}$
deformations (2-cells) in $\c$, for $0\leq i<j<k\leq p$.

The other two simplicial sets in the diagram, $\gner\c:\Delta^{\!{^\mathrm{op}}} \to \set$ and
$\nabla\c:\Delta^{\!{^\mathrm{op}}} \to \set$, here termed  `geometric nerves' of $\c$, are
respective non-normalized versions of the above unitary ones by Street and Duskin. And the
remaining five `nerves' in the diagram associate to every bicategory $\c$ simplicial objects in
$\cat$: Those called its `unitary categorical geometric nerves', which are those denoted by
$\cgner\hspace{-4pt}^{^\mathrm{u}}\c$ and $\overline{\nabla}\hspace{-2.5pt}_{^{_\mathrm{u}}}\c$,
its `categorical geometric nerves', denoted by $\cgner\c$ and $\overline{\nabla}\c$, and its `Segal
nerve', written as $\underline{\mathrm{S}}\c$, respectively. As we will prove, all these are
`special' simplicial categories, in the sense  that \cite{segal74, thomason} the Segal projection
maps on them induce homotopy equivalences on classifying spaces. Furthermore, the latter,
$\underline{\mathrm{S}}\c$, is a weak 2-category in the sense of Tamsamani \cite{tam} and Simpson
\cite{sim}, that is, the Segal projection maps on it are surjective equivalences of categories as
observed by Lack and Paoli in \cite{lack-paoli} (where $\underline{\mathrm{S}}\c$ is called the
`2-nerve of $\c$', but this may be  confusing terminology since, for example, the unitary geometric
nerve  $\cgner\hspace{-4pt}^{^\mathrm{u}}\c$ of a 2-category $\c$ is also called the `2-nerve of
$\c$' in  \cite{tonks}).

Our first major result with this work states and proves that: \begin{quote}{\em For any bicategory
$\c$, all the continuous maps in the diagram
$$
\xymatrix{&&\ \  \class\c \ar@<1.5pt>[d]\ar@<1.5pt>[dl]\ \ \ar@<-1.5pt>[dr]&&\\ \class\cgner\c\ &
\ar@<1.5pt>[l]\ \class\cgner\hspace{-3pt}^{^\mathrm{u}}\c\ &  \ar@<1.5pt>[l]\
\class\underline{\mathrm{S}}\c \ \ar@<-1.5pt>[r] &\
\class\overline{\nabla}\hspace{-2.5pt}_{^{_\mathrm{u}}}\c\  \ar@<-1.5pt>[r]& \
\class\overline{\nabla}\c
\\ \ar@<1.5pt>[u]\class\gner\c\
&\ar@<1.5pt>[l]\ \class\gner\hspace{-3pt}^{^\mathrm{u}}\c \ \ar@<1.5pt>[u] &&\
\class\nabla\hspace{-2.5pt}_{^{_\mathrm{u}}}\c \ \ar@<1.5pt>[u] \ar@<-1.5pt>[r] & \ \class\nabla\c,
\ \ar@<1.5pt>[u]}
$$
induced by $(A)$ on classifying spaces, are homotopy equivalences.} \end{quote}

Many properties of the classifying space construction for bicategories, $\c\mapsto \class\c$, may
be easier to establish depending on the nerve used for realizations. For example, in what is our
second main result, we deal with the homotopy theory of diagrams of bicategories and homomorphisms
$\c:I^{^{\,\mathrm{op}}}\to \bicat$. Following Grothendieck \cite{grothendieck}, we show the
existence of a `bicategorical Grothendieck construction' $\int_I \!\c, $
 that suitably assembles all bicategories $\c_i$, $i\in\mathrm{Ob}I$, and whose classifying  space $\class\!\int_I \!\c$ can be thought as the homotopy colimit of the diagram of spaces $i\mapsto \class\c_i$, as we prove by using geometric nerves for realizing classifying spaces of bicategories. More precisely, we prove the following:
\begin{quote}
{\em Suppose a category $I$ is given. For every functor $\c:I^{^{\mathrm{op}}}\to\bicat$, there
exists a natural weak homotopy equivalence of simplicial sets
$$ \xymatrix{\eta: \hoco_I\gner\c\longrightarrow
\gner\!\int_I\!\c,}$$  where $\gner\c: I^{^{\,\mathrm{op}}}\to\sset$ is the diagram of simplicial
sets obtained by composing $\c$ with the geometric nerve of bicategories functor. }
\end{quote}

We should stress that, regarding any category as a discrete bicategory, the above weak equivalence
$\eta$  for diagrams of categories $\c:I^{^{\,\mathrm{op}}}\to \cat$
  just gives the weak equivalence shown  by Thomason  in his well-known Homotopy Colimit Theorem \cite{thomason}.  However, our proof that $\eta$ is a weak homotopy equivalence follows different lines than Thomason's.

The plan of this paper is, briefly, as follows. After this introductory Section 1, the paper is
organized in six sections. Section 2 aims to make this paper as self contained as possible; hence,
at the same time as we fix notations and terminology, we also review some necessary aspects and
results from the background of simplicial sets, homotopy theory of categories and bicategories.
However, the material in this Section 2 is quite standard, so the expert reader may skip most of
it. In Section 3, we use Jardines's supercoherence theorem \cite{jardine}  to introduce the
pseudo-simplicial nerve of a bicategory $\cner\c$ and then  Thomason's theory on lax-diagrams of
categories define its classifying space $\class \c$. The following two sections 4 and 5 are
dedicated to describing  the remaining nerves of a bicategory discussed in the paper as well as the
simplicial maps, functors, or pseudo-functors that connect to each other. In Section 6, we include
our main result on the homotopy invariance of all the different nerves of a bicategory, and the
final Section 7 mainly contains our homotopy colimit theorem for diagrams of bicategories.

\section{Preliminaries}

In this preliminary section we review, without any claim to originality, some standard
constructions and well-known basic facts mainly concerning nerves and classifying spaces of (small)
categories. In Subsection 2.4, we fix some terminology and notational conventions on bicategories
that for some readers may be idiosyncratic.

\subsection{Some facts concerning  simplicial sets}
Hereafter, we shall regard each ordered set $[n]=\{0,1,\ldots,n\}$ as the category with exactly one
arrow $j\rightarrow i$ if $i\leq j$. Then, a non-decreasing map $[n]\rightarrow [m]$ is the same as
a functor, so that we see $\Delta $, the simplicial category of finite ordinal numbers, as a full
subcategory of $\cat$, the category (actually the 2-category) of small categories.

The category of {\em simplicial sets}, that is, the category of functors
${S:\Delta^{\!{^\mathrm{op}}}\to \set}$, where $\set$ is the category of sets, is denoted by
$\sset$.  The simplicial standard $n$-simplex
$\Delta[n]=\mbox{Hom}_{\Delta^{\!{^\mathrm{op}}}}(-,[n])$ is the obvious representable functor, as
usual.

Recall that the category $\Delta$ is generated by the injections $d^i:[n-1]\to[n]$ (cofaces),
$0\leq i\leq n$, which omit  the $i$th element and the surjections $s^i:[n+1]\to [n]$
(codegeneracies), $0\leq i\leq n$, which repeat the $i$th element, subjet to the well-known  {\em
cosimplicial identities}: $d^j d^i=d^i d^{j-1}$ if $i<j$, etc. Thus, in order to define a {\em
simplicial object} in a category $\mathcal{E}$, say $S:\Delta^{\!{^\mathrm{op}}}\to  \mathcal{E}$,
it suffices to give the objects (of $n$-simplices) $S_n$, $n\geq 0$, together with morphisms
$$
\begin{array}{lll}d_i=(d^i)^*=S(d^i):S_n\to S_{n-1},&\ 0\leq i\leq n& \ (\mbox{the face operators}),\\[3pt]
s_i=(s^i)^*=S(s^i):S_n\to S_{n+1},&\ 0\leq i\leq n& \ (\mbox{the degeneracy operators}),
\end{array}
$$
satisfying the well-known basic {\em simplicial identities}: $d_id_j=d_{j-1}d_i$ if $i<j$, etc. A
{\em simplicial morphism} $f:S\rightarrow S^\prime$ is just a natural transformation from $S$ to
$S^\prime$; it then consists of a family $\{f_n:S_n\rightarrow S^\prime_{n},\, n\geq 0\}$ of arrows
in $\mathcal{E}$ that commute with the face and degeneracy operators. If $f,g:S\rightarrow
S^\prime$ are simplicial morphisms, then a \emph{simplicial homotopy} $H:f\Rightarrow g$ is  a
system $\{H_m:S_n\rightarrow S_{n+1}^\prime,\; 0\leq m\leq n\}$ of arrows in $\mathcal{E}$ that
satisfies the set of homotopy identities described, for example, in \cite[Definition 5.1]{may67}.

Segal's {\em geometric realization} \cite{segal} of a simplicial (compactly generated topological)
space $S:\Delta^{\!{^\mathrm{op}}}\to  \mathbf{Top}$ is denoted by $\class S$. This construction is
functorial and any simplicial homotopy $H:f\Rightarrow g$, between simplicial space maps
$f,g:S\rightarrow S^\prime$, determines a homotopy $\class H:\class f \Rightarrow \class g$,
\cite[Corollary 11.10]{may72}. For instance, by regarding a set as a discrete space, the (Milnor's)
geometric realization of a simplicial set $S:\Delta^{\!{^\mathrm{op}}}\to \set$ is $\class S$. If
$f,g:S\to S'$ are simplicial maps, between simplicial sets,  then a simplicial homotopy
$H:f\Rightarrow g$ amounts to the same thing \cite[Proposition 6.2]{may67} as a simplicial map $H:
S\times \Delta[1]\to S'$ making this diagram commutative:
$$
\xymatrix@C=0pt{S\times\Delta[0]\ar[d]_{1\times d^1}&\cong S \ar[rrrd]^{f}&\\ S\times\Delta[1] \ar[rrrr]^H&&&&S'.\\
S\times\Delta[0]\ar[u]^{1\times d^0}&\cong S \ar[urrr]_{g}& }
$$

A {\em weak homotopy equivalence} of simplicial sets is a simplicial map whose geometric
realization is a homotopy equivalence or, equivalently, induces isomorphisms in homotopy groups.

A \emph{bisimplicial set} is a functor
$S:\Delta^{\!{^\mathrm{op}}}\times\Delta^{\!{^\mathrm{op}}}\rightarrow\set$. This amounts to a
family of sets $\{S_{p,q};\,  p,q\geq 0\}$ together with horizontal and vertical face and
degeneracy operators
$$\xymatrix{S_{p+1,q}&\ar[l]_-{\textstyle s_i^h}S_{p,q}\ar[r]^-{\textstyle d_i^h}&S_{p-1,q}, \hspace{0.6cm}
S_{p,q+1}&\ar[l]_-{\textstyle s_j^v} S_{p,q}\ar[r]^-{\textstyle d_j^v}&S_{p,q-1},}
$$
with $0\leq i\leq p$ and $0\leq j\leq q$ respectively, such that, for all $p$ and $q$, both
$S_{p,\ast}$ and $S_{\ast,q}$ are simplicial sets and the horizontal operators commute with the
vertical ones. Note that, on the one hand, any bisimplicial set $S$ provides two simplicial objects
in the category of simplicial sets: the horizontal one $S^h:[p]\mapsto S_{p,\ast}$ and the vertical
one $S^v:[q]\mapsto S_{\ast,q}$. Then, by taking realization, $S$ gives rise to two simplicial
spaces $\class S^h:[p]\mapsto \class S_{p,\ast}$ and $\class S^v:[q]\mapsto \class S_{\ast,q}$,
respectively. On the other hand, by composing with the diagonal functor
$\diag:\Delta^{\!{^\mathrm{op}}}\rightarrow\Delta^{\!{^\mathrm{op}}}\times\Delta^{\!{^\mathrm{op}}}$,
the bisimplicial set $S$ also provides another simplicial set $\diag S:[n]\mapsto S_{n,n}$, whose
face and degeneracy operators are given in terms of those of $S$ by the formulas $d_i=d_i^hd_i^v$
and $s_i=s_i^hs_i^v$, respectively. It is known (e.g. \cite[Lemma in page 86]{quillen}) that there
are natural homeomorphisms
\begin{equation}\label{homeo}
\class\class S^h \cong \class\diag S \cong \class\class S^v.
\end{equation}

It is a relevant fact  that,  if $f:S\rightarrow S^\prime$ is a bisimplicial map such that the
induced maps $ f_{p,\ast}: S_{p,\ast}\rightarrow  S^\prime_{p,\ast}$ (respect. $ f_{\ast,q}:
S_{\ast,q}\rightarrow S^\prime_{\ast,q}$) are weak homotopy equivalences for all $p$ (respect.
$q$), then so is the map $\diag f: \diag S \rightarrow \diag S^\prime$ (see \cite[Chapter XII, 4.2
and 4.3]{bousfield-kan} or \cite[IV, Proposition 1.7]{g-j}, for example).

We shall also use the bar construction on a bisimplicial set $\overline{W}S$, also called its {\em
codiagonal} or {\em total complex}. Let us recall that the functor
$$
\overline{W}: \bsset \to \sset
$$
is  the right Kan extension along the ordinal sum functor $\gner\times \gner \to\gner$,
$([p],[q])\mapsto [p+1+q]$. For any given bisimplicial set $S$, $\overline{W}S$ can be described as
follows \cite[\S III]{artin-mazur}: the set of $p$-simplices of $\overline{W}S$ is
$$
\Big\{(t_{0,p} \dots,t_{p,0})\in
\prod_{m=0}^{p}S_{m,p-m}~|~d^v_0t_{m,p-m}=d^h_{m+1}t_{m+1,p-m-1},\, 0\leq m< p\Big\}$$ and, for
$0\leq i\leq p$, the faces and degeneracies of a $p$-simplex are given by
$$
\begin{array}{l}d_i(t_{0,p}
\dots,t_{p,0})=(d_i^vt_{0,p},\dots,d_i^vt_{i-1,p-i+1},d_i^ht_{i+1,p-i-1},\dots,d_i^ht_{p,0}),
\\[-4pt]~\\
s_i(t_{0,p} \dots,t_{p,0})=(s_i^vt_{0,p},\dots,s_0^vt_{i,p-i},s_i^ht_{i,p-i},\dots,s_i^ht_{p,0}).
\end{array}
$$

For any bisimplicial set $S$, there is a natural weak homotopy equivalence \cite{cegarra1,cegarra2}
\begin{equation}\label{Phi} \Phi:\diag\,S\to \overline{W}S,
\end{equation}
which carries a $p$-simplex $t_{p,p}\in \diag \,S$ to
$$
\Phi t_{p,p}=\Big((d_1^h)^pt_{p,p},(d_2^h)^{p-1}d_0^vt_{p,p},\dots,
(d_{m+1}^h)^{p-m}(d_0^v)^mt_{p,p},\dots, (d_0^v)^pt_{p,p}\Big).
$$

\subsection{The classifying space of a small category}
In Quillen's development of $K$-theory, the higher $K$-groups are defined as the homotopy groups of
a topological {\em classifying space}, $\class \c$, functorially associated to a small category
$\c$, the geometric realization of its so-called nerve (Grothendieck, \cite{grothendieck}). This
nerve is a simplicial set encoding the structure of the category in terms of its faces and
degeneracies, and it can be easily described by means of two isomorphic constructions:

On one hand, the nerve of a category $\c$ can be defined as the simplicial set
\begin{equation}\label{p.1.1}
\ner\c :\Delta^{\!{^\mathrm{op}}}\longrightarrow  \set, \end{equation} whose set of $p$-simplices
$$\begin{array}{l} \ner\c_p = \bigsqcup\limits _{(x_0,\ldots,x_p)\in \mbox{\scriptsize
Ob}\c^{p+1}}\hspace{-0.6cm}
\c(x_1,x_0)\times\c(x_2,x_1)\times\cdots\times\c(x_p,x_{p-1}),\\[15pt] \ner\c_0=\mbox{Ob}\c,\end{array}$$ consists of  length $p$ sequences of composible morphisms in $\c$
$$\xymatrix@C=15pt{x_0&\ar[l]_<<{\textstyle u_1}x_1&\ar[l]_<<{\textstyle u_2}\cdots&\ar[l]_<<{\textstyle u_p}x_p}.
$$
The face and degeneracy operators are defined by the well-known formulas:
\begin{equation}\label{p.1.2}\begin{array}{l} d_i(u_1,\dots,u_p)=\left\{
\begin{array}{lcl}
\hspace{-0.2cm}  (u_2,\dots,u_p) &\text{ if }& i=0 ,  \\[4pt]
\hspace{-0.2cm}  (u_1,\dots,u_i\,u_{i\text{+}1},\dots,u_p) & \text{ if } & 0<i<p , \\[4pt]
\hspace{-0.2cm}  (u_0,\dots,u_{p-1}) & \text{ if }& i=p,
\end{array}\right.\\[20pt]
s_i(u_1,\dots,u_p)=(u_1,\dots,u_i,1_{x_i},u_{i\text{+}1},\dots, u_p).
\end{array}
\end{equation}

On the other hand, let
\begin{equation}\label{p.1.3} \begin{array}{ll}\Delta\c:&\Delta^{\!{^\mathrm{op}}}\longrightarrow Set,\\[6pt] &[p]\mapsto \func([p],\c),\end{array} \end{equation} be the simplicial set whose $p$ simplices are the functors $F:[p]\to {\mathcal C}$ or, equivalently, tuples of
arrows in $\c$ $$F=(\xymatrix{Fj\ar[r]^{\textstyle F_{i,j}}& Fi})_{_{0\leq i\leq j\leq p}}$$ such
that $F_{i,j}F_{j,k}=F_{ik}$, for $i\leq j\leq k$ and $F_{i,i}=1_{Fi}$.

As there is quite an obvious simplicial isomorphism
$$ \begin{array}{l}\Delta\c\cong \ner\c,\\[6pt] F\mapsto
(\xymatrix@C=25pt{F0&\ar[l]_{\textstyle F_{0,1}}F1&\ar[l]_{\textstyle F_{1,2}}\cdots&\ar[l]_{
\textstyle F_{p-1,p}}Fp})
\end{array} $$
both simplicial sets $\Delta\c$ and $\ner\c$ are usually identified, and hereafter we will also do
so. However, note that $\Delta\c$ has a more pleasing geometric interpretation than $\ner\c$, since
a $p$-simplex of $\Delta\c$ can be thought of as the 1-skeleton of an oriented standard $p$-simplex
with objects $Fi$ of $\c$ placed on the vertices and arrows $F_{i,j}:Fj\to Fi$ placed on the edges
for $0\leq i < j\leq p$, and the requirement that every triangle
$$\xymatrix@C=7pt@R=7pt{ & Fj\ar[dr]^{\textstyle F_{i,j}} &\\
Fk \ar[rr]_{\textstyle F_{i,k}} \ar[ur]^{\textstyle F_{j,k}} && Fi }$$ be commutative.

For instance, note that $\ner[n]\cong \Delta[n]$ is the simplicial standard $n$-simplex, and
$\class[n]=\Delta_n$ is the standard $n$-dimensional affine simplex. When a monoid (or group)
${\mathcal M}$ is regarded as a category with only one object, then $\class{\mathcal M}$ is its
classifying space in the traditional sense. Therefore, many weak homotopy types thus occur, since
every path-connected space has the weak homotopy type  of the classifying space of a monoid
(\cite[Theorem 3.5]{fiedorowicz}). Moreover, any CW-complex is homotopy equivalent to the
classifying space of a small category, as Quillen showed \cite{quillen}: The category of simplices
$\int_\Delta {\hspace{-2pt}}S$, of a simplicial set $S$, has as objects pairs $(p,x)$ where $p\geq
0$ and $x$ is a $p$-simplex of $S$; and arrow $\alpha:(p,x)\rightarrow (q,y)$ is an arrow
$\alpha:[p]\to [q]$ in $\Delta$ with the property $x=\alpha^*y$. Then there exists a homotopy
equivalence $|S|\simeq \class\int_\Delta {\hspace{-2pt}}S$ between the geometric realization of $S$
and the classifying space of $\int_\Delta{\hspace{-2pt}}S$ (this result is, in fact, a very
particular case of the homotopy colimit theorem of Thomason \cite{thomason}). If $X$ is any
CW-complex and we take $S=SX$, the total singular complex of $X$, then $X\simeq |SX|\simeq
B\int_\Delta{\hspace{-2pt}} SX$.

The nerve construction on categories gives a fully faithful functor, embedding the category $\cat$
into the category of simplicial sets. Moreover, this functor obviously commutes with products and,
as pointed out by Segal \cite[Proposition 2.1]{segal}, these facts imply that, for any two functors
$F,G:{\mathcal C}\to {\mathcal D}$, the functor nerve defines a bijection
$$\text{Nat. transformations} (F,G)\cong ~\text{Simpl. homotopies}(\Delta F, \Delta G),$$
between the set of natural transformations $\alpha:F\Rightarrow G$ and the set of simplicial
homotopies between the induced simplicial maps on the nerves: Interpret $\alpha$ as a functor
$\alpha:{\mathcal C} \times [1]\rightarrow {\mathcal D}$. Hence, what $\alpha$ induces is a
simplicial map $\Delta \alpha:\Delta{\mathcal C}\times \Delta[1]\to \Delta{\mathcal D}$, that is, a
simplicial homotopy. As a consequence, functors related by some natural transformation go to
homotopic cellular maps on classifying spaces. In particular, if a functor $F:{\mathcal C}\to
{\mathcal D}$ has a left or right adjoint, the induced map $\class F:\class{\mathcal C}\to \class
{\mathcal D}$ is a homotopy equivalence.

\subsection{The classifying space of a diagram of  categories}
In \cite{segal}, Segal extended Milnor's geometric realization process to simplicial (compactly
generated topological) spaces and provided, for instance, the notion of classifying spaces for
simplicial categories $\c:\Delta^{\!{^\mathrm{op}}}\to \cat$: By replacing each category $\c_p$,
$p\geq 0$, by its classifying space $\class\c_p$, one obtains a simplicial space, whose Segal
realization is, by definition, the classifying space $\class\c$ of the simplicial category $\c$,
that is,
\begin{equation}\label{p.1.8} \class\c=\class([p]\mapsto \class\c_p). \end{equation} But
note, as an instance of the homeomorphisms (\ref{homeo}), that there is a natural homeomorphism
\begin{equation}\label{p.1.9} \class\c \cong \class\diag \Delta \c \end{equation}
 with the geometric realization of the simplicial set diagonal of the bisimplicial set
obtained by composing $\c:\Delta^{\!{^\mathrm{op}}}\to \cat$ with the nerve functor ${\Delta:\cat \to\sset} $, that is, $$\begin{array}{ll}\Delta\c:&\Delta^{\!{^\mathrm{op}}}\times \Delta^{\!{^\mathrm{op}}}\longrightarrow \set.\\[4pt]
&([p],[q])\longmapsto \func([q],\c_p)\end{array}$$

Segal's construction above, for simplicial categories, is actually  a particular case of the more
general notion of classifying space for arbitrary {\em diagrams of categories}: If $\c:
I^{{^\mathrm{op}}}\to \cat$ is a functor,  where $I$ is any category, then its classifying space
\begin{equation}\label{p.1.11} \class\c=\class\,\hoco_I\Delta\c, \end{equation}
is given through the homotopy colimit construction by Bousfield and Kan \cite[\S
XII]{bousfield-kan}, of the diagram of simplicial sets $\Delta \c:I^{{^\mathrm{op}}}\to \sset$,
obtained by composing ${\c}$ with the nerve functor. That is, the simplicial set
$$ \begin{array}{ll} \hoco_I\Delta{\c}:&   \Delta^{\!{^\mathrm{op}}}\longrightarrow \set. \\[4pt]
& [p] \mapsto \bigsqcup\limits_{[p]\overset{\sigma}\to I} \func([p],\c_{\sigma 0})\end{array} $$

Since, for any simplicial category $\c:\Delta^{\!{^\mathrm{op}}}\to \cat$, there is a natural weak
homotopy equivalence of simplicial sets \cite[XII, 4.3]{bousfield-kan},
\begin{equation}\label{bk}\mbox{hocolim}_\Delta\Delta\c\overset{\thicksim}\longrightarrow \mbox{diag}\Delta{\mathcal C},\end{equation} then
both constructions (\ref{p.1.8}) and (\ref{p.1.11}) for the classifying space $\class\c$ of a
simplicial category ${\mathcal C}$ coincide up to a natural homotopy equivalence.

\subsection{Some bicategorical conventions}

 We employ the standard nomenclature concerning bicategories, and we refer to \cite{benabou,duskin,g-p-s,gurski} or \cite{street} for the background. For the sake of clarity, we state the following:

 A small {\em bicategory}, $\c$, provides us with the following data: (i)  a set $\mbox{Ob}\c$ of {\em objects} (or 0-{\em cells}) of $\c$; (ii) for each ordered pair of objects $(y,x)$ of $\c$, a category $\c(y,x)$ whose objects $u:y\to x$ are called {\em morphisms} (or 1-{\em cells}) of $\c$ with source $y$ and target $x$, and whose arrows ${\alpha:u\Rightarrow u'}$  are called {\em deformations} (or 2-{\em cells}) of $\c$ and are usually depicted as $\xymatrix @C=0pc {y  \ar@/^0.5pc/[rr]^u \ar@/_0.5pc/[rr]_{u'} & {}_{\Downarrow\alpha} &x }$. The composition of deformations in each category $\c(y,x)$ is called the {\em vertical composition} and will be denoted by $\beta\cdot \alpha$; (iii) for every triplet $(z,y,x)$ of objects of $\c$, a functor $\c(y,x)\times \c(z,y)
\overset{\circ}\to \c(z,x)$,
$$\xymatrix @C=0pc{x&\   {}_{\textstyle \Downarrow\alpha} & \ar@/^0.8pc/[ll]^{\textstyle u'} \ar@/_0.8pc/[ll]_{\textstyle u} y&\ {}_{\textstyle \Downarrow\alpha'}&\ar@/^0.8pc/[ll]^{\textstyle v'} \ar@/_0.8pc/[ll]_{\textstyle v}  z}\overset{\textstyle \circ}\mapsto \xymatrix@C=0pc{x& \  {}_{\textstyle \Downarrow{\alpha\circ \alpha'}} &\ar@/^0.9pc/[ll]^{\textstyle  u'\circ v'} \ar@/_0.9pc/[ll]_{\textstyle u\circ v} z },$$ called {\em horizontal composition}; (iv) for every object $x$ of $\c$, a distinguished morphism $1_x:x\to x$, called the {\em identity} of $x$; (v) invertible deformations  $\boldsymbol{a}_{u,v,w}:u\circ(v\circ w)\overset{\thicksim}\Rightarrow (u\circ v)\circ w$, called the {\em associativity constraints}, which are natural in $(u,v,w)\in \c(y,x)\times\c(z,y)\times\c(t,z)$; and (vi) invertible deformations $\boldsymbol{l}_u: 1_x\circ u \overset{\thicksim}\Rightarrow u$ and
$\boldsymbol{r}_u:u\circ 1_y\overset{\thicksim}\Rightarrow u$, called {\em identity constraints},
which are natural in $u\in \c(y,x)$. These constraints are required to be coherent, in the sense
that the following pentagons and triangles are commutative:
\begin{equation}\label{pent}
\xymatrix@C=18pt{u_1\circ( u_2\circ (u_3\circ u_4)) \ar@2{->}[r]^{\textstyle
\boldsymbol{a}}\ar@2{->}[d]_{\textstyle 1\circ \boldsymbol{a}}&(u_1\circ u_2)\circ (u_3\circ u_4
)\ar@2{->}[r]^{\textstyle \boldsymbol{a}}& ((u_1\circ u_2)\circ u_3)\circ u_4
\\ u_1\circ ((u_2\circ u_3)\circ u_4)\ar@2{->}[rr]^{\textstyle \boldsymbol{a}}&&(u_1\circ(u_2\circ u_3))\circ u_4\,,\ar@2{->}[u]_{\textstyle \boldsymbol{a}\circ 1}}\end{equation}
\begin{equation}\label{tri1}\xymatrix@C=10pt{u\circ (1_y\circ v)\ar@2{->}[rr]^{\textstyle \boldsymbol{a}}\ar@2{->}[rd]_{\textstyle 1 \circ \boldsymbol{l}}&&(u\circ 1_y)\circ v
\ar@2{->}[ld]^{\textstyle \boldsymbol{r}\circ 1}\\&u\circ v\,.&}
\end{equation}
In any bicategory, for any object $x$, the equality
\begin{equation}\label{rl}\boldsymbol{l}_{1_x}=\boldsymbol{r}_{1_x}\end{equation}
holds, and, for any two morphisms $z\overset{v}\to y\overset{u}\to x$, the two triangles below
commute.
\begin{equation}\label{tri23}\xymatrix@C=10pt{u\circ (v\circ 1_z)\ar@2{->}[rr]^{\textstyle \boldsymbol{a}}\ar@2{->}[rd]_{\textstyle \boldsymbol{1\circ r}}&&(u\circ v)\circ 1_z
\ar@2{->}[ld]^{\textstyle \boldsymbol{r}}\\&u\circ v&}\hspace{0.5cm} \xymatrix@C=10pt{1_x\circ (u
\circ v)\ar@2{->}[rr]^{\textstyle \boldsymbol{a}}\ar@2{->}[rd]_{\textstyle \boldsymbol{l}}&&(1_x
\circ u)\circ v \ar@2{->}[ld]^{\textstyle \boldsymbol{l}\circ 1}\\&u\circ v&}
\end{equation}

A bicategory in which all the constraints are identities is a 2-{\em category}, that is, just a
category enriched in the category $\cat$ of small categories. As each category $\c$ can be
considered as a 2-category in which all deformations are identities, that is, in which each
category $\c(x,y)$ is discrete, several times throughout the paper categories are considered as
special bicategories.

If $\b$, $\c$ are  bicategories, then a {\em lax functor} $ F=(F,\widehat{F}):\b \rightsquigarrow
\c$ consists of: (i) a mapping $F:\mbox{Ob}\b\to \mbox{Ob}\c$; (ii)
 for each ordered pair of objects objects $(y,x)$ of $\b$, a functor $F:\b(y,x)\to \c(Fy,Fx)$; (iii) deformations $\widehat{F}_{u,v}:Fu\circ Fv\Rightarrow F(u\circ v)$ that are natural in $(u,v)\in\b(y,x)\times\b(z,y)$; and (iv) for each object $x$ of $\b$, a deformation $\widehat{F}_x:1_{Fx}\Rightarrow F1_x$. All these data are subject to the coherence commutativity conditions:
$$
\xymatrix{ Fu\circ (Fv\circ Fw)\ar@{=>}[r]^{\textstyle 1\circ \widehat{F}}\ar@{=>}[d]^{\textstyle
\boldsymbol{a}}&Fu\circ F(v\circ w)\ar@{=>}[r]^{\textstyle \widehat{F}}&
F(u\circ(v\circ w))\ar@{=>}[d]_{\textstyle F\boldsymbol{a}}\\
(Fu\circ Fv)\circ Fw\ar@{=>}[r]^{\textstyle \widehat{F}\circ 1}&F(u\circ v)\circ Fw
\ar@{=>}[r]^{\textstyle \widehat{F}}&F((u\circ v)\circ w),}
 $$
 $$
 \xymatrix{Fu\circ 1_{Fy}\ar@{=>}[r]^{\textstyle \boldsymbol{r}}\ar@{=>}[d]_{\textstyle 1\circ \widehat{F}}&Fu
 \\ Fu\circ F1_y\ar@{=>}[r]^{\textstyle \widehat{F}}&F(u\circ 1_y),\ar@{=>}[u]_{\textstyle F\boldsymbol{r}}}\hspace{0.5cm}
 \xymatrix{1_{Fx}\circ Fu\ar@{=>}[r]^{\textstyle \boldsymbol{l}}\ar@{=>}[d]_{\textstyle \widehat{F}\circ 1}&Fu
 \\ F1_x\circ Fu\ar@{=>}[r]^{\textstyle \widehat{F}}&F(1_x\circ u).\ar@{=>}[u]_{\textstyle F\boldsymbol{l}}}
 $$

 Replacing the structure deformations above with $\widehat{F}_{u,v}:F(u\circ v)\Rightarrow Fu\circ Fv$ and $\widehat{F}_x: F(1_x)\Rightarrow 1_{Fx}$, we have the notion of {\em oplax functor} $F: \b\rightsquigarrow \c$.

A lax functor is termed a {\em pseudo functor} or a {\em homomorphism} whenever all the structure
constraints $\widehat{F}_{u,v}:Fu\circ Fv\Rightarrow F(u\circ v)$ and
$\widehat{F}_x:1_{Fx}\Rightarrow F(1_x)$ are invertible. When these deformations are all
identities, then $F$ is called a 2-{\em functor} and is written as $F:\b\to\c$. If the unit
constraints $\widehat{F}_x$ are all identities, then the lax functor is qualified as (strictly)
{\em unitary} or {\em normal}.

The composition of lax functors $\a \overset{F}\rightsquigarrow \b\overset{G}\rightsquigarrow \c$
will be denoted by juxtaposition, that is, $\a \overset{GF}\rightsquigarrow \c$. Recall that its
constraints are obtained from those of $F$ and $G$ by the rule $\widehat{GF}=G\widehat{F}\cdot
\widehat{G}F$; that is, by the compositions
$$
\begin{array}{l}
\widehat{GF}_{u,v}: \xymatrix{GFu\circ GFv\ar@{=>}[r]^{\textstyle \widehat{G}_{Fu,Fv}}&G(Fu\circ Fv)\ar@{=>}[r]^{\textstyle G\widehat{F}_{u,v}}&GF(u\circ v)},\\[5pt]
\widehat{GF}_{x}:\xymatrix{1_{GFx}\ar@{=>}[r]^{\textstyle
\widehat{G}_{Fx}}&G1_{Fx}\ar@{=>}[r]^{\textstyle G\widehat{F}_x}&GF1_x}.
\end{array}
$$
The composition of lax functors is associative and unitary, so that the category of bicategories
and lax functors is defined. The subcategory of bicategories with homomorphisms between them will
be denoted by $\bicat$.

If $F,F':\b \rightsquigarrow \c$ are lax functors, then we follow the convention of \cite{g-p-s} in
what is meant by a {\em lax transformation} ${\alpha=(\alpha,\widehat{\alpha}):F\Rightarrow F'}$.
Thus, $\alpha$ consists of morphisms ${\alpha x:Fx\to F'x}$, $x\in \mbox{Ob}\b$, and deformations
$$
\xymatrix@C=7pt@R=7pt{&Fx\ar[rd]^{\textstyle \alpha x}&\\Fy\ar[ru]^{\textstyle
Fu}\ar[rd]_{\textstyle \alpha y}&~~ \Downarrow\widehat{\alpha}_u&F'\!x\\&F'\!y\ar[ru]_{\textstyle
F'\!u}&}
$$
that are natural on morphisms $u:y\to x$, subject to the usual two commutativity axioms:
$$
\xymatrix{\alpha x\circ(Fu\circ Fv)\ar@{=>}[r]^{\textstyle 1\circ \widehat{F}} \ar@{=>}[d]_{\textstyle \boldsymbol{a}} &\alpha x\circ F(u\circ v)\ar@{=>}[r]^{\textstyle \widehat{\alpha}_{u\circ v}}&F'(u\circ v)\circ \alpha z\\
(\alpha x\circ Fu)\circ Fv \ar@{=>}[d]_{\textstyle \hat{\alpha}_u\circ 1} &&(F'\!u\circ F'\!v)\circ \alpha z \ar@{=>}[u]_{\textstyle \widehat{F'}\circ 1}\\
(F'\!u\circ \alpha y)\circ Fv \ar@{=>}[r]^{\textstyle \boldsymbol{a}^{-1}}&F'\!u\circ(\alpha y\circ
Fv)\ar@{=>}[r]^{\textstyle  1\circ \widehat{\alpha}_v} & F'\!u\circ (F'v\circ \alpha z)\,,
\ar@{=>}[u]_{\textstyle \boldsymbol{a}}}
$$
$$
\xymatrix{\alpha x\circ F1_x\ar@{=>}[rr]^{\textstyle \widehat{\alpha}_{1_x}}& &F'\!1_x\circ \alpha
x\\\alpha x\circ 1_{Fx}\ar@{=>}[u]^{\textstyle 1\circ \widehat{F}} \ar@{=>}[r]^{\textstyle
\boldsymbol{r}}&\alpha x\ar@{=>}[r]^{\textstyle \boldsymbol{l}^{-1}}&1_{F'\!x}\circ \alpha x\,.
\ar@{=>}[u]_{\textstyle \widehat{F'}\circ 1}}
$$

Replacing the structure deformation above with $\widehat{\alpha}_u: F'\!u\circ \alpha y\Rightarrow
\alpha x\circ Fu$, we have the notion of {\em oplax transformation} $\alpha:F\Rightarrow F'$.

\subsection{The classifying space of a lax diagram of categories}
In nature, actual functors ${\mathcal C}:I^{^{\mathrm{op}}}\to \cat$ are rare, but pseudo-functors
are ubiquitous. The above construction (\ref{p.1.11}), of classifying spaces for diagrams of
categories, was extended to lax diagrams of categories by Thomason \cite{thomason} using methods by
Grothendieck. Recall that, for any given small category $I$,  a {\em lax diagram of categories}
means a lax functor
$$\c=(\c,\widehat{\c}):I^{^{\,\mathrm{op}}} \rightsquigarrow\cat,$$  to the 2-category $\cat$ of small categories, functors, and natural transformations. So $\c$ is a system of data consisting of a category $\c_i$ for each object $i$ of $I$, a functor
$\c_i\overset{a^*}\to \c_j$ for each arrow $j\overset{a}\to i$ of $I$, a natural transformation
$\widehat{\c}=\widehat{\c}_{a,b}:b^*a^*\Rightarrow (ab)^*$, for each two composible arrows
$k\overset{b}\to j\overset{a}\to i$ in  $I$, and a natural transformation
$\widehat{\c}=\widehat{\c}_i:1_{\c_i}\Rightarrow 1_i^*$, for each object $i$ of $I$. These must
satisfy the conditions that, for $\ell \overset{c}\to k\overset{b}\to j \overset{a}\to i$, the
following diagram commutes:
$$\xymatrix@C=7pt@R=9pt{&c^*b^*a^*\ar@2{->}[dl]_{\textstyle c^*\widehat{\c}}
\ar@2{->}[dr]^{\textstyle \widehat{\c} a^*
}&\\
c^*(ab)^*\ar@2{->}[dr]_{\textstyle ~\widehat{\c}}& &(bc)^*a^*
\ar@2{->}[dl]^{\textstyle ~\widehat{\c}},\\
&(abc)^*&}
$$
and, for $j\overset{a}\to i$, the compositions $\xymatrix@C=15pt{a^*\ar@{=>}[r]^{\textstyle
\widehat{\c} a^*}&1_j^*a^*\ar@{=>}[r]^{\textstyle \widehat{\c}}&a^*}$ and
$\xymatrix@C=15pt{a^*\ar@{=>}[r]^{\textstyle a^*\widehat{\c}}&a^*1_i^*\ar@{=>}[r]^{\textstyle
\widehat{\c}}&a^*}$ are both the identity transformation on the functor $a^*:\c_i\to\c_j$.

 The so-called {\em Grothendieck
construction} on a lax diagram of categories ${\c:I^{^{\mathrm{op}}} \rightsquigarrow\cat}$,
denoted by $$ \xymatrix{\int_I\c,} $$ is the category whose objects are pairs $(x,i)$, where $i$ is
an object of $I$ and $x$ is one of ${\mathcal C}_i$; a morphism $(u,a):(y,j)\to (x,i)$ in
$\int_I\c$ is a pair of morphisms where $a:j\to i$ in $I$ and $u:y\to a^*x$ in $\c_j$. If
$(v,b):(z,k)\to (y,j)$ is another morphism in $\int_I\c$, then we have the morphisms
$$\xymatrix{z\ar[r]^-{\textstyle v}& b^*y\ar[r]^-{\textstyle b^*u}&
b^*a^*x\ar[r]^{\textstyle \widehat{\c}x}&(ab)^*x}$$ in $\c_k$, and the composition of $(v,b)$ with
$(u,a)$ is defined by
$$ (u,a)(v,b) =(\widehat{\c}x\cdot  b^*\!u\cdot  v, ab):(z,k)\to (x,i).$$
The identity morphism of an object $(x,i)$ is $(x\overset{\textstyle
\widehat{\c}_ix}\longrightarrow 1_i^*x,1_i)$.

Then, the classifying space of the lax diagram $\c:I^{{^\mathrm{op}}} \rightsquigarrow\cat$,
$\class \c$,  is defined to be the classifying space of its Grothendieck construction, that is,
\begin{equation}\label{p.1.15}\xymatrix{\class\c=\class \int_I\c}.\end{equation}

 Thomason's homotopy colimit theorem \cite[Theorem 1.2]{thomason}
states the following:

\begin{theorem}\label{thoma} Let $\c:I^{{^\mathrm{op}}} \to \cat$ be a functor. There is a natural weak homotopy equivalence
\begin{equation}\label{p.1.18}\xymatrix{ \eta:\hoco_I\ner\c \longrightarrow \ner\int_I\c.} \end{equation}
\end{theorem}

 Therefore, both definitions (\ref{p.1.11}) and (\ref{p.1.15}), for the classifying space $\class\c $ of a diagram of categories $\c:I^{{^\mathrm{op}}} \to \cat$, lead to the
same space, up to a natural homotopy equivalence. Furthermore, the construction of $\class\c$ is
consistent with the so-called {\em Street rectification} process, $\c \mapsto \c'$ \cite{street72}
(see also \cite[Theorem 3.4]{may}), which associates to any lax diagram of categories
$\c:I^{^{\mathrm{op}}}\rightsquigarrow\cat$ a homotopy equivalent strict diagram
$\c':I^{^{\mathrm{op}}}\to \cat$, that is, such that there is a natural homotopy equivalence
\cite[Lemma 3.2.5]{thomason} $$ \class\c\simeq\class \c'. $$

The Grothendieck construction on lax diagrams $\c:I^{^{\mathrm{op}}}\rightsquigarrow\cat$ is
natural both in $I$ and $\c$. Recall that, given lax diagrams $\c,\d:
I^{^{\mathrm{op}}}\rightsquigarrow\cat$, a {\em lax morphism} (or oplax transformation)
$$F=(F,\widehat{F}):\c\rightsquigarrow\d,$$
is a system of data consisting of a functor $F_i:\c_i\to \d_i, $ for each object $i$ of $I$, and a
natural transformation $\widehat{F}_a:F_j a^*\Rightarrow a^*F_i$,
$$
\xymatrix@C=8pt@R=8pt{&\c_j\ar[rd]^{\textstyle F_j}& \\ \c_i\ar[rd]_{\textstyle
F_i}\ar[ru]^{\textstyle a^*}&\Downarrow \!\widehat{F}_a &\d_j,\\ &\d_i\ar[ru]_{\textstyle a^*}&}
$$
for each morphism $j\overset{a}\to i$ in $I$. These are subject to the conditions that the
following diagrams commute:
$$
\xymatrix{F_k b^*a^*\ar@{=>}[r]^{\textstyle \widehat{F}_b a^*}\ar@{=>}[d]_{\textstyle F_k\widehat{\c}}&b^*F_j a^*\ar@{=>}[r]^{\textstyle b^*\widehat{F}_a}&b^*a^*F_i\ar@{=>}[d]^{\textstyle \widehat{\c} F_i}\\
F_k(ab)^*\ar@{=>}[rr]^{\textstyle \widehat{F}_{ab}}&&(ab)^*F_i,}\
\xymatrix{&F_i1_i^*\ar@{=>}[rd]^{\textstyle \widehat{F}_{1_i}}&\\F_i\ar@{=>}[rr]^{\textstyle
\widehat{\c} F_i}\ar@{=>}[ru]^{\textstyle F_i\widehat{\c}}&&1_i^*F_i,}
$$
for every pair of composible morphisms $k\overset{b}\to j \overset{a}\to i$ and any object $i$ in
$I$.  Any lax morphism $F:\c\rightsquigarrow \d$ induces a functor $\int_IF:\int_I\c\to \int_I\d$,
whence a cellular map $$\class F:\class \c\to \class \d,$$ defined on objects by
$F_*(x,i)=(F_ix,i)$. For a morphism $(u,a):(y,j)\to (x,i)$ in $\int_I\c$, we have the composible
morphisms in $\d_j$
$$ F_jy\overset{\textstyle F_ju}\longrightarrow F_j a^*x\overset{\textstyle \widehat{F}_a x}\longrightarrow a^*F_ix,$$
and $F_*(u,a)=(\widehat{F}_a x \cdot F_ju,a)$. A main result by Thomason \cite[Corollary
3.3.1]{thomason} states the following:
\begin{theorem}\label{p.1.23} If $F:\c\rightsquigarrow\d$ is a lax morphism between lax diagrams $\c,\d:  I^{^{\mathrm{op}}}\rightsquigarrow\cat$ such that the induced maps $BF_i:\class \c_i\to \class\d_i$ are homotopy equivalences, for all objects $i$ of
$I$, then the induced map $\class F:\class\c \to \class\d$ is a homotopy equivalence.
\end{theorem}

\section{The pseudo-simplicial nerve and the classifying space of a bicategory.}

Let $\c$ be any given bicategory. When $\c$ is strict (i.e., a 2-category),  then the nerve
construction (\ref{p.1.1}) actually works by giving a simplicial category $\cner
\c:\Delta^{\!^{\mathrm{op}}}\to \cat$, whose Segal's classifying space (\ref{p.1.8}), or
(\ref{p.1.9}), is usually taken to be the classifying space of the 2-category (see \cite{b-c,
hinich, moerdijk-svensson, til, til2} or \cite{thomason}, for examples). For an arbitrary $\c$, its
classifying space is defined in a similar way as for the strict case; however, the process is more
complicated since the horizontal composition in a bicategory is in general not associative and not
unitary (which is crucial for constructing the simplicial category $\cner\c$) but it is only so up
to coherent isomorphisms. This  `defect'  has the effect of forcing one to deal with the
classifying space of a nerve of $\c$, which is not simplicial but only up to coherent isomorphisms:
the (normal) pseudo-simplicial category

\begin{equation}\label{ps1.1}\cner\c: \Delta^{\!^{\mathrm{op}}}\rightsquigarrow \cat
\end{equation} whose category of $p$-simplices is
$$
\cner\c_p = \bigsqcup_{(x_0,\ldots,x_p)\in \mbox{\scriptsize Ob}\c^{p+1}}\hspace{-0.3cm}
\c(x_1,x_0)\times\c(x_2,x_1)\times\cdots\times\c(x_p,x_{p-1}),
$$
where a typical arrow is a string of deformations in $\c$
\begin{equation}\label{xi}
\xi=\xymatrix @C=1pc {x_0  & {\Downarrow\, \alpha_1} & x_1 \ar@/^1pc/[ll]^{\textstyle v_1}
\ar@/_1pc/[ll]_-{\textstyle u_1} & {\Downarrow\, \alpha_2} & x_2  \ar@/^1pc/[ll]^{\textstyle v_2}
\ar@/_1pc/[ll]_-{\textstyle u_2}&\hspace{-15pt}\cdots&\hspace{-10pt}
 x_{n-1}
  & {\Downarrow\, \alpha_n} & x_n
\ar@/^1pc/[ll]^{\textstyle v_n} \ar@/_1pc/[ll]_-{\textstyle u_n}},
\end{equation}
and $\cner\c_0=\mbox{Ob}\,\c$, as a discrete category.

The face and degeneracy functors are defined in the standard way by the formulas (\ref{p.1.2}),
both for morphisms and deformations, substituting  juxtaposition with the symbol $\circ$, used for
the horizontal composition in $\c$. That is, $d_0(\xi)= (\alpha_2,\dots,\alpha_n)$,
$d_1(\xi)=(\alpha_1\circ \alpha_2,\dots,\alpha_n)$, and so on.

If $a:[q]\to[p]$ is any non-identity map in $\Delta$, then we write $a$ in the (unique) form (see
\cite{may67}, for example) $a=d^{i_1}\cdots d^{i_s}s^{j_1}\cdots s^{j_t}$, where $0\leq i_s<\cdots<
i_1 \leq p$, $0\leq j_1<\cdots< j_t \leq q$ and $q+s=p+t$, and the induced functor
$a^*:\cner_p\c\to\cner_q\c$ is defined by $a^*=s_{j_t}\cdots s_{j_1}d_{i_s}\cdots d_{i_1}$. Note
that $d_jd_i=d_id_{j+1}$ for $i\leq j$, unless $i=j$ and $1\leq i\leq p-2$, in which case the
associativity constraint of $\c$ gives a  canonical natural isomorphism
\begin{equation}\label{ps1.2} d_id_i \cong d_id_{i+1}. \end{equation}
Similarly, all the equalities $d_0s_0=1$, $d_{p+1}s_p=1$, $d_is_j=s_{j-1}d_i$ if $i<j$ and
$d_is_j=s_jd_{i-1}$ if $i>j+1$, hold, and the unit constraints of $\c$ give canonical isomorphisms
\begin{equation}\label{ps1.3} d_is_i\cong 1, \hspace{0.6cm} d_is_{i+1}\cong 1. \end{equation}

Then it is a fact that this family of natural isomorphisms (\ref{ps1.2}) and (\ref{ps1.3}) uniquely
determines a whole system of natural isomorphisms $$ b^*a^*\cong (ab)^*,$$ one for each pair of
composible maps in $\Delta$, $[n]\overset{b}\to [q]\overset{a}\to [p]$, such that the assignments
$\alpha\mapsto \alpha^*$, $1_{[p]}\mapsto 1_{\cner\c_p}$, together with these isomorphisms
$b^*a^*\cong (ab)^*$, give the data for the pseudo-simplicial category (\ref{ps1.1}),
$\cner\c:\Delta^{\!^{\mathrm{op}}}\rightsquigarrow\cat$. This fact can be easily proven by using
Jardine's supercoherence theorem \cite[Corollary 1.6]{jardine}, since the commutativity of the
seventeen diagrams of supercoherence, (1.4.1)-(1.4.17) in \cite{jardine},
 easily follows from the coherence conditions (\ref{pent}),  (\ref{tri1}),  (\ref{rl}), and  (\ref{tri23}).

Recalling the construction (\ref{p.1.15}) for the classifying spaces of lax diagrams of categories,
we state the following:

\begin{definition}\label{ps1.5}
The classifying space $\class\c$, of a bicategory $\c$, is the classifying space of its
pseudo-simplicial nerve {\em (\ref{ps1.1})}, $\cner\c:\Delta^{\!^{\mathrm{op}}}\rightsquigarrow
\cat$, that is,
\begin{equation}\label{ps1.6} \xymatrix{\class\c=\class\int_\Delta \cner\c}.
\end{equation}
\end{definition}

\begin{remark}\label{2-cat}{\em Let $\c$ be a 2-category. Then, its pseudo-simplicial nerve (\ref{ps1.1}) is actually a
simplicial category $$\cner\c:\Delta^{\!^{\mathrm{op}}}\to\cat,$$ and there are natural homotopy
equivalences
$$\class\c\overset{(\ref{p.1.18})}\simeq |\hoco_\Delta\cner\c|\overset{(\ref{bk})}\simeq|\diag\ner\cner\c|\overset{(\ref{p.1.9})}\simeq \class\cner\c,$$
where $\ner\cner\c: ([p],[q])\mapsto \ner(\cner\c_p)_q$ is the bisimplicial set {\em double nerve}
of the 2-category obtained by composing $\cner\c$ with the functor nerve of categories, and
$\class\cner\c$ is the classifying space (\ref{p.1.8}) of the simplicial category $\cner\c$, that
is,  the Segal realization  of the simplicial space $[p]\mapsto \class\cner\c_p$.\qed }
\end{remark}
\begin{remark}\label{mono}{\em A {\em monoidal} (tensor) category $\m=(\m,\otimes,\text{I},\boldsymbol{a},\boldsymbol{l},\boldsymbol{r})$, \cite{maclane},  can be viewed as a
bicategory $$\Omega^{^{-1}}\hspace{-5pt}\m$$
 with only one object, say $*$,  the objects $u$ of $\m$ as morphisms $u:*\rightarrow *$ and the morphisms of $\m$
as deformations. Thus, $\Omega^{^{-1}}\hspace{-5pt}\m(*,*)=\m$, and it is the horizontal
composition of morphisms and deformations given by the  tensor functor
${\otimes:\m\times\m\rightarrow\m}$. The identity at the object is $1_*=\text{I}$, the unit object
of the monoidal category, and the associativity, left unit, and right unit constraints for
$\Omega^{^{-1}}\hspace{-5pt}\m$ are just those of the monoidal category, that is,
$\boldsymbol{a}$, $\boldsymbol{l}$ and $\boldsymbol{r}$, respectively.

 The pseudo-simplicial nerve of the bicategory $\Omega^{^{-1}}\hspace{-5pt}\m$, as in (\ref{ps1.1}), is
 exactly the pseudo-simplicial category that the monoidal category  defines by the reduced bar construction; that is, the pseudo-simplicial category
$$\cner\Omega^{^{-1}}\hspace{-5pt}\m:\Delta^{\!^{\mathrm{op}}}
\rightsquigarrow\cat,\hspace{0.5cm}[p]\mapsto \m^p,$$ whose category of $p$-simplices is the
$p$-fold power of the underlying category $\m$, with faces and degeneracy functors defined by
analogy with those of the nerve of a monoid. Therefore, {\em the classifying space  of the monoidal
category} (see \cite[\S 3]{jardine}, \cite[Appendix]{hinich} or \cite{b-c2}, for example) is just
$$\class \Omega^{^{-1}}\hspace{-5pt}\m,$$ the classifying space of the one object bicategory it
defines.

The observation, due to Benabou \cite{benabou}, that monoidal categories are essentially the same
as bicategories with just one object, is known as the {\em delooping principle}, and the bicategory
$\Omega^{^{-1}}\hspace{-5pt}\m$ is called the {\em delooping of the monoidal category}
\cite[2.10]{k-v}. The reason for this terminology  is the existence of a natural map $$\class \m\to
\Omega(\class\Omega^{^{-1}}\hspace{-5pt}\m,*),$$  where $\class\m$ is the classifying space of the
underlying category $\m$, which is, up to group completion, a homotopy equivalence (see
\cite[Propositions 3.5 and 3.8]{jardine} or \cite[Corollary 4]{b-c2}, for example; also, see Remark
\ref{gp} for a proof therein). Then, the higher $K$-groups $K_i$, $i>0$, of the monoidal category
$\m$ are the $(i+1)$-th homotopy groups of its classifying space $\class
\Omega^{^{-1}}\hspace{-5pt}\m$ (cf.\cite{segal74}). When the monoidal category has a given
symmetry, then $\class \Omega^{^{-1}}\hspace{-5pt}\m$ is precisely the space at level 1 of the
$\Omega$-spectrum associated to the symmetric monoidal category (see \cite[4.2.2]{thomason} or
\cite[3.12]{jardine}, for example). The most striking instance is $\m=A$, the strict monoidal
category with only one object defined by an abelian group $A$, where both compositions and tensor
products are given by the addition in $A$; in this case, $\class A$ is a $K(A,1)$-space, and
$\Omega^{^{-1}}\hspace{-5pt} A$ is a bicategory with only one object and only one arrow whose
classifying space $\class\Omega^{^{-1}}\hspace{-5pt}A$ is a $K(A,2)$-space. \qed}
\end{remark}

Any homomorphism $F:\b\rightsquigarrow\c$, between bicategories, gives rise to a morphism of
supercoherent structures in the sense of Jardine \cite{jardine}, $F_*:\cner\b\to \cner\c$, that, on
a morphism $\xi$ as in (\ref{xi}), of the category of $p$-simplices of $\cner\b$, acts by
$$\xymatrix@C=0.5pt@R=0.5pt{
\xi \ar@{|-{>}}[rrrr]^{F_*} &&~&& Fx_0  &{\Downarrow F\alpha_1} & Fx_1\ar@/^1pc/[ll]^{\textstyle
Fv_1} \ar@/_1pc/[ll]_-{\textstyle Fu_1}  & {\Downarrow F\alpha_2} & Fx_2 \ar@/^1pc/[ll]^{\textstyle
Fv_2} \ar@/_1pc/[ll]_-{\textstyle Fu_2}\ar@{.}[rr] && Fx_{p-1}  & {\Downarrow F\alpha_p} &
Fx_p.\ar@/^1pc/[ll]^{\textstyle Fv_p} \ar@/_1pc/[ll]_-{\textstyle Fu_p}}$$ The structure natural
isomorphisms $s_i F_* \cong F_*s_i$ and $d_iF_*\cong F_* d_i$ are canonically obtained from  the
invertible  structure constraints of the homomorphism, $\widehat{F}:1_{Fx_i}\cong F1_{x_i}$ and
$\widehat{F}: F(u_i)\circ F(u_{i+1})\cong F(u_i\circ u_{i+1})$ (the commutativity of the needed six
coherence diagrams in \cite{jardine} is clear). Then, $F:\b\rightsquigarrow\c$ determines a
pseudo-simplicial functor, $\cner F:\cner\b\rightsquigarrow \cner\c$, and therefore a functor
$\int_\Delta \!F:\int_\Delta\!\b\to \int_\Delta\!\c$ and a corresponding map on the classifying
spaces $$ \class F:\class\b\to\class\c. $$

Thus, the classifying space construction (\ref{ps1.6}), $\c\to\class \c$, defines a functor from
bicategories, with homomorphisms between them, to CW-complexes. However, as we will see later, any
lax, or oplax, functor $F:\b\rightsquigarrow\c$, that is, without the requirement that its
structure constraints $\widehat{F}$ be isomorphisms, also induces a continuous map $\class
F:\class\b\to\class\c$, well defined up to homotopy equivalence.

\section{The simplicial sets geometric nerves of a bicategory.}
The two isomorphic constructions (\ref{p.1.1}) and (\ref{p.1.3}) to define the nerve, and then the
classifying space, of a category suggest different extensions to bicategories. Looking at
(\ref{p.1.1}), we were led to definition (\ref{ps1.1}), which recovers the more traditional way
through which a classifying space is assigned in the literature to certain kinds of bicategories,
such as 2-categories (see Remark \ref{2-cat}) or monoidal categories (see Remark \ref{mono}). But
the construction we have taken, in Definition  \ref{ps1.5}, for the classifying space of a
bicategory $\class \c$ runs through the Grothendieck construction on its pseudo-simplicial nerve,
which implies that the cells of $\class\c$ have little apparent intuitive connection with the cells
of the original bicategory and that they do not enjoy any proper geometric meaning. However,
looking at (\ref{p.1.3}), we are led to dealing with another convincing way of associating a space
to any bicategory $\c$: through its {\em unitary geometric nerve}  as defined by Street in
\cite{street} (see also \cite{duskin} and \cite{gurski2}).
\begin{definition} The unitary geometric nerve of a bicategory $\c$ is the simplicial set
\begin{equation}\label{ngn} \begin{array}{ll}\Delta^{\hspace{-2pt}^\mathrm{u}}\!\c:&\Delta^{\!^{\mathrm{op}}}\ \to \ \set,\\[6pt]
&[p]\mapsto \mathrm{NorLaxFunc}([p],\c)\end{array}\end{equation} whose $p$-simplices are the normal
lax functors ${F:[p]\rightsquigarrow \c}$. If $a:[q]\to [p]$ is any map in $\Delta$, that is, a
functor,  the induced
$a^*:\Delta^{\hspace{-2pt}^\mathrm{u}}\!\c_p\to\Delta^{\hspace{-2pt}^\mathrm{u}}\!\c_q$ carries
${F:[p]\rightsquigarrow \c}$ to the composite ${Fa:[q]\rightsquigarrow \c}$, of $F$ with $a$.
\end{definition}

This simplicial set $\Delta^{\hspace{-2pt}^\mathrm{u}}\c$ encodes the entire bicategorical
structure of $\c$, and the following lemma allows us to show a pleasing geometrical description of
its simplices:

\begin{lemma}\label{gn.1} Let $\c$ be a bicategory. Any system of data consisting of objects $Fi$, $0\leq i\leq p$, morphisms  $F_{i,j}:Fj\to Fi$, $0\leq i<j\leq p$,  and deformations ${\widehat{F}_{i,j,k}:F_{i,j}\circ F_{j,k}\Rightarrow F_{i,k}}$, $0\leq i<j<k\leq p$, such that, for ${0\leq i<j<k<\ell\leq p}$, the  following square of deformations
$$\xymatrix@C=50pt{F_{i,j}\circ ( F_{j,k}\circ F_{k,\ell})\ar@{=>}[rr]^{\textstyle 1\circ \widehat{F}_{{j,k,\ell}}}
\ar@{=>}[d]_{\textstyle \boldsymbol{a}}&& F_{i,j}\circ F_{j,\ell} \ar@{=>}[d]^{\textstyle \widehat{F}_{{i,j,\ell}}} \\
(F_{i,j}\circ F_{j,k})\circ F_{k,\ell}\ar@{=>}[r]^{\textstyle\widehat{F}_{{i,j,k}}\circ 1}&
F_{i,k}\circ F_{k,\ell} \ar@{=>}[r]^-{\textstyle \widehat{F}_{{i,k,\ell}}}&F_{i,\ell}
 }
 $$
commutes in the category $\c(F\ell,Fi)$,  uniquely extends to a normal lax functor
${F:[p]\rightsquigarrow \c}$.
\end{lemma}
\begin{proof} The whole data for $F:[p]\rightsquigarrow \c$ are obtained by putting $F_{i,i}=1_{Fi}$, $\widehat{F}_{i,i,j}=\boldsymbol{l}:1_{Fi}\circ F_{i,j}\Rightarrow F_{i,j}$, ${\widehat{F}_{i,j,j}=\boldsymbol{r}:F_{i,j}\circ 1_{Fj}\Rightarrow F_{i,j}}$, and $\widehat{F}_i=1_{1_{Fi}}:1_{Fi}\Rightarrow 1_{Fi}$. So defined, $F$ is actually a (normal) lax functor thanks to the commutativity of the coherence triangles in (\ref{tri1}) and (\ref{tri23}), and the equality in (\ref{rl}). The uniqueness of $F$ is clear.
\end{proof}

Thus, for a bicategory  $\c$, the vertices of its normal geometric nerve
$\Delta^{\hspace{-2pt}^\mathrm{u}}\c$ are  the objects $F0$ of $\c$, the 1-simplices are the
morphisms $\xymatrix{F0&\ar[l]_-{\textstyle F_{0,1}}F1}$,  and the 2-simplices are triangles
$$
\xymatrix{\ar@{}[drr]|(.6)*+{\Downarrow {\textstyle \widehat{F}_{_{0,1,2}}}}  & F1 \ar[dl]_{\textstyle F_{0,1}}            \\
 F0  & &  \ar[ll]^{\textstyle F_{0,2}}    F2\ar[ul]_{\textstyle F_{1,2}}, }
$$
with $\widehat{F}_{0,1,2}:F_{0,1}\circ F_{1,2}\Rightarrow F_{0,2}$ a deformation in $\c$. For
$p\geq 3$, a $p$-simplex of $\Delta^{\hspace{-2pt}^\mathrm{u}}\c$ is geometrically represented by a
diagram in $\c$ with the shape of the 2-skeleton of an oriented standard $p$-simplex, whose faces
are triangles
$$
\xymatrix{\ar@{}[drr]|(.6)*+{\Downarrow \widehat{F}_{_{i,j,k}}}
                & Fj \ar[dl]_{\textstyle F_{i,j}}             \\
Fi  & &     Fk   \ar[ul]_{\textstyle F_{j,k}} \ar[ll]^{\textstyle F_{i,k}}     }
$$
with objects $Fi$ placed on the vertices, $0\leq i\leq p$, morphisms $F_{i,j}:Fj\rightarrow Fi$ on
the edges, $0\leq i<j\leq p$, and $\widehat{F}_{i,j,k}:F_{i,j}\circ F_{j,k}\Rightarrow F_{i,k}$
deformations, for $0\leq i<j<k\leq p$. These data are required to satisfy the condition that each
tetrahedron
$$
\begin{array}{ccc}
\xymatrix {
 & F\ell \ar[dl]_{\textstyle F_{i,\ell}}  \ar[dd]|<<<<<<{\textstyle F_{j,\ell} }\ar[dr]^{\textstyle F_{k,\ell}}& \\
Fi   & & Fk \ar[dl]^{\textstyle F_{j,k}} \ar[ll]|<<<<<<{\; \textstyle F_{i,k}\;}\\
 & Fj \ar[ul]^{\textstyle
F_{i,j}}}
 & \hspace{1cm}&
\xymatrix@R=2pt{\\  \widehat{F}_{i,j,k}:F_{i,j}\circ F_{j,k}\Rightarrow F_{i,k}\\
\widehat{F}_{i,j,\ell}:F_{i,j}\circ F_{j,\ell}\Rightarrow F_{i,\ell}\\
\widehat{F}_{i,k,\ell}:F_{i,k}\circ F_{k,\ell}\Rightarrow F_{i,\ell}\\
\widehat{F}_{j,k,\ell}:F_{j,k}\circ F_{k,\ell}\Rightarrow F_{j,\ell},}
\end{array}
$$
for $0\leq i< j< k< \ell \leq p$ is commutative in the sense that the following equation on
deformations holds:
$$
\xymatrix{&F\ell\ar[dd]\ar[rd]\ar[ld] &&&& F\ell\ar[rd]\ar[ld]\ar@{}[d]|(0.6)*+{
\overset{\widehat{F}_{_{i,k,\ell}}}\Leftarrow}&\\ Fi &\ar@{}[l]|(.4)*+{
\overset{\widehat{F}_{_{i,j,\ell}}}\Leftarrow}&Fk\ar[ld]
\ar@{}[l]|(.6)*+{ \overset{\widehat{F}_{_{j,k,\ell}}}\Leftarrow}&=&Fi&&Fk.\ar[ld]\ar[ll]\\
&Fj\ar[lu] &&&&Fj\ar[lu]\ar@{}[u]|(.6)*+{ \overset{\widehat{F}_{_{i,j,k}}}\Uparrow}&}
$$

The simplicial set $\Delta^{\hspace{-2pt}^\mathrm{u}}\c$ becomes coskeletal in dimensions greater
than 3. More precisely, for $p\geq 3$, a $p$-simplex $F:[p]\rightsquigarrow \c$ of
$\Delta^{\hspace{-2pt}^\mathrm{u}}\c$ is determined uniquely by its boundary $\partial
F=(d_0F,\dots,d_pF)$
 $$\xymatrix{\partial\Delta[p]\ar[r]^{\partial F}\ar@{_{(}->}[d]&\Delta^{\hspace{-2pt}^\mathrm{u}}\c,\\
\Delta[p]\ar[ru]_F&}$$
 and,  for $p\geq 4$, every possible boundary of a $p$-simplex, $\partial\Delta[p]\to \Delta^{\hspace{-2pt}^\mathrm{u}}\c$, is actually the boundary $\partial F$ of a $p$-simplex $F$.

 For several discussions, it is suitable to handle the (non-normalized) {\em geometric nerve} of a bicategory $\c$:
\begin{definition} The geometric nerve of a bicategory $\c$ is the simplicial set
\begin{equation}\label{gn} \begin{array}{ll}\Delta\c:&\Delta^{\!^{\mathrm{op}}}\ \to \ \set,\\[6pt]
&[p]\mapsto \lfunc([p],\c),\end{array}\end{equation} that is, the simplicial set whose
$p$-simplices are all lax functors ${F:[p]\rightsquigarrow \c}$.\end{definition}
  Hence, the unitary geometric nerve  $\Delta^{\hspace{-2pt}^\mathrm{u}}\c$ becomes a simplicial subset of $\Delta\c$.  The $p$-simplices of the  geometric nerve $\Delta\c$ are described similarly  to those of the normalized one, but now they include  deformations of
$\c$ \begin{equation}\label{1.2.16'} \widehat{F}_i:1_{F_i}\Rightarrow F_{i,i},\hspace{0.6cm}0\leq
i\leq p, \end{equation} with the requirement that the diagrams
$$
\xymatrix@C=9pt@R=9pt{ F_{i,j}\circ F_{j,j}\ar@2{->}[dr]_{\textstyle \widehat{F}_{i,j,j}} &&
F_{i,j}\circ
1_{F_j}\ar@2{->}[dl]^{\textstyle \boldsymbol{r}}\ar@2{->}[ll]_{\textstyle 1\circ \widehat{F}_j}\\
& F_{i,j}&\\}  ~~~~ \xymatrix@C=9pt@R=9pt{ F_{i,i}\circ F_{i,j}\ar@2{->}[dr]_{\textstyle
\widehat{F}_{i,i,j}} &&\ar@2{->}[ll]_{\textstyle \widehat{F}_i\circ 1} 1_{F_i}\circ
F_{i,j}\ar@2{->}[dl]^{\textstyle \boldsymbol{l}}\\
& F_{i,j}&\\}
$$
be commutative.

Note   that the geometric nerve construction on bicategories, $\c\mapsto \Delta \c$, is clearly
functorial on lax functors between bicategories,  whereas $\Delta^{\hspace{-2pt}^\mathrm{u}}\c$ is
functorial only on normalized lax functors between bicategories.

We have used the above lax functors $[p]\rightsquigarrow \c$ to define geometric nerves of
bicategories. However, there is no substantial reason for not considering  oplax functors. In fact,
we also have {\em geometric nerves}
\begin{equation}\label{ogn} \begin{array}{ll}\nabla\c:&\Delta^{\!^{\mathrm{op}}}\ \to \ \set,\\[6pt]
&[p]\mapsto {\ensuremath{\mathrm{OpLaxFunc}}}([p],\c)\end{array}\end{equation} that is, the
simplicial set whose $p$-simplices are all oplax functors, ${F:[p]\rightsquigarrow \c}$, and
\begin{equation}\label{ongn} \begin{array}{ll}\nabla\hspace{-2.5pt}_{^{_\mathrm{u}}}\c:&\Delta^{\!^{\mathrm{op}}}\ \to \ \set,\\[6pt]
&[p]\mapsto \mathrm{NorOpLaxFunc}([p],\c)\end{array}\end{equation}
 the simplicial subset of $\nabla\c$,  whose $p$-simplices are the normal oplax functors from $[p]$ to $\c$.

\begin{remark}{\em In \cite{duskin}, Duskin gave  a characterization of the normal geometric nerve $\nabla\hspace{-2.5pt}_{^{_\mathrm{u}}}\c$ of a bicategory $\c$ in terms of its simplicial structure. The result states that a simplicial set is  isomorphic to the normal geometric nerve of a bicategory if and only if it satisfies the coskeletal conditions above as well as supporting appropriate sets of `abstractly invertible' 1- and 2-simplices. For  instance, recall that a {\em bigroupoid} is a bicategory $\c$ in which every deformation is invertible, that is, all categories $\c(y,x)$ are groupoids, and every morphism $u:y\to x$ is a biequivalence, that is,  there exist  a morphism $u':x\to y$ and deformations $u\circ u'\Rightarrow 1_y$ and $1_x\Rightarrow u'\circ u$. Then, the normal geometric nerve of a bicategory $\c$ becomes a Kan complex \cite{may} (fibrant, \cite{g-j}), that is, every extension problem $$\xymatrix{\Lambda^k[n]\ar[r]\ar@{_{(}->}[d]&\nabla\hspace{-2.5pt}_{^{_\mathrm{u}}}\c\\
\Delta[n]\ar@{..>}[ru]&}$$ of any $(k,n)$-horn, $\Lambda^k[n]\to
\nabla\hspace{-2.5pt}_{^{_\mathrm{u}}}\c$,  has a solution if and only if the bicategory $\c$ is a
bigroupoid. In such a case, the homotopy groups of $\nabla\hspace{-2.5pt}_{^{_\mathrm{u}}}\c$ are

\vspace{0.2cm}
$\pi_i(\nabla\hspace{-2.5pt}_{^{_\mathrm{u}}}\c,x)$ $$=\left\{\begin{array}{lll}K_0\c&\text{the set of iso-classes of objects of } \c& \text{ for }i=0,\\
 K_0\c(x,x)&\text{the group of iso-classes of automorphisms in $\c$ of } x&\text{ for } i=1,\\K_1\c(x,x)&\text{the abelian group of autodeformations in $\c$ of } 1_x &\text{ for }i=2,\\
 0&\text{ for }i\geq 3.\end{array}\right.$$

Furthermore, the normal geometric nerve construction embeds the category of bigroupoids, with
normal homomorphisms between them, as a full and reflexive subcategory of the category of Kan
complexes, whose replete image consists of those Kan complexes $K$ for which every $(k,n)$-horn
$\Lambda^k[n]\to K$ has a unique filler
$$\xymatrix{\Lambda^k[n]\ar[r]\ar@{_{(}->}[d]&K,\\
\Delta[n]\ar@{..>}[ru]_{\exists !}&}$$ for all $0\leq k\leq n$ and $n\geq 3$; that is, $K$ is a Kan
2-{\em hypergropoid} in the sense of Duskin-Glenn. }\qed
\end{remark}

\begin{remark}{\em
The normal geometric nerve of a 2-category (even of an $n$-category) was first studied in
\cite{street2}.  In \cite{tonks}, but under the name of `2-nerve of 2-categories', the normal
geometric nerve construction $\Delta^{\hspace{-2pt}^\mathrm{u}}$ was considered for proving that
the category $\mathbf{2}\text{-}\cat$, of small 2-categories and 2-functors, has a Quillen model
structure such that the functor $Ex^2\,\Delta^{\hspace{-2pt}^\mathrm{u}}:\mathbf{2}\text{-}\cat\to
\sset$ is a right Quillen equivalence of model categories, where $Ex$ is the endofunctor in $\sset$
right adjoint to the subdivision $Sd$ (see \cite{g-j}, for example). In \cite{b-c},  it was proved
that, for any 2-category $\c$, there is a natural homotopy equivalence $\class\c\simeq
\class\Delta^{\hspace{-2pt}^\mathrm{u}}\c$, and therefore it follows that the correspondence
$\c\mapsto\class\c$  induces an equivalence between the corresponding homotopy category of
2-categories and the ordinary homotopy category of CW-complexes. By this correspondence,
2-groupoids correspond to spaces whose homotopy groups $\pi_n$ are trivial for $n > 2$
\cite{moerdijk-svensson}, and from this point of view the use of 2-categories and their classifying
spaces in homotopy theory goes back to Whitehead (1949) and Mac Lane-Whitehead (1950), since
2-groupoids with only one object (2-groups) are the same as crossed modules. }\qed
\end{remark}

\begin{remark}{\em \label{e2} Let $\m$ be a monoidal category. The delooping bicategory  $\Omega^{^{-1}}\hspace{-5pt}\m$ (see Remark  \ref{mono}) realizes the classifying space of the monoidal category and, in \cite{b-c2}, it is shown how that realization can be made through the normal geometric nerve construction. That is, there is a homotopy equivalence $\class \Omega^{^{-1}}\hspace{-5pt}\m\simeq \class\nabla\hspace{-2.5pt}_{^{_\mathrm{u}}}\Omega^{^{-1}}\hspace{-5pt}\m$.

This geometric nerve $\nabla\hspace{-2.5pt}_{^{_\mathrm{u}}}\Omega^{^{-1}}\hspace{-5pt}\m$ is a
$3$-coskeletal reduced (1-vertex) simplicial set, whose simplices have the following simplified
interpretation: the $1$-simplices  are the objects $F_{0,1}$ of $\m$, the $2$-simplices are
morphisms of $\m$ of the form
$$
\xymatrix@C=30pt{ F_{0,2}\ar[r]^-{\textstyle F_{0,1,2}}  &  F_{0,1}\otimes F_{1,2} }
$$
and the $3$-simplices are commutative diagrams in $\m$ of the form
$$
\xymatrix@C=40pt{F_{0,3}\ar@{=>}[r]^-{\textstyle {F}_{{0,2,3}}} \ar@{=>}[d]_{\textstyle
{F}_{{0,1,3}}}&F_{0,2}\otimes F_{2,3} \ar@{=>}[r]^-{\textstyle {F}_{{0,1,2}}\otimes 1}
&(F_{0,1}\otimes F_{1,2})\otimes F_{2,3}\\ F_{0,2}\otimes F_{1,3}\ar@{=>}[rr]^{\textstyle 1\otimes
F_{{1,2,3}}}& &F_{0,1}\otimes (F_{1,2}\otimes F_{2,3}). \ar@{=>}[u]_{\textstyle \boldsymbol{a}} }
$$

 Bigroupoids with only one object are delooping bicategories  $\Omega^{^{-1}}\hspace{-5pt}\m$ of {\em categorical groups} \cite{joyal} (Gr-categories, \cite{bre92}), that is, monoidal groupoids $\m$ in which translations are autoequivalences. Normal geometric nerves of categorical groups are studied in \cite{Ca-Ceg, ceg-gar}.
 }\qed
\end{remark}

\section{More $\cat$-valued nerves of a bicategory.}

In this section, we shall describe the following commutative diagram with the various `nervesï of a
bicategory $\c$ discussed in the paper:

\begin{equation}\label{dia}
\xymatrix{&&\ \  \cner_{_\text{}}\c \ar@{_{(}~>}@<1.5pt>_{J}[d]\ar@{_{(}~>}@<1.5pt>[dl]_J\ \
\ar@{^{(}~>}@<-1.5pt>[dr]^J&&\\ \cgner\c\ & \ar@{_{(}->}@<1.5pt>[l]\
\cgner\hspace{-4pt}^{^\mathrm{u}}\c_{_\text{}}\ &  \ar@{_{(}->}@<1.5pt>[l]\
\underline{\mathrm{S}}\c \ \ar@{^{(}->}@<-1.5pt>[r] &\
\overline{\nabla}\hspace{-2.5pt}_{^{_\mathrm{u}}}\c\  \ar@{^{(}->}@<-1.5pt>[r]& \
\overline{\nabla}\c
\\ \ar@{_{(}->}@<1.5pt>[u]\gner\c\
&\ar@{_{(}->}@<1.5pt>[l]\ \gner\hspace{-4pt}^{^\mathrm{u}}\c \ \ar@{_{(}->}@<1.5pt>[u] &&\
\nabla\hspace{-2.5pt}_{^{_\mathrm{u}}}\c \ \ar@{_{(}->}@<1.5pt>[u] \ar@{^{(}->}@<-1.5pt>[r] & \
\nabla\c, \ \ar@{_{(}->}@<1.5pt>[u]}
\end{equation}
where $\underline{\ner}\c$, in the top row, is the pseudo simplicial nerve (\ref{ps1.1}) of the
bicategory, and $\gner\c$, $\gner\hspace{-4pt}^{^\mathrm{u}}\!\c$, $\nabla\c$ and
$\nabla\hspace{-2.5pt}_{^{_\mathrm{u}}}\c$, in the bottom row, are the geometric nerves  of the
bicategory, as defined in (\ref{gn}), (\ref{ngn}), (\ref{ogn}), and (\ref{ongn}) respectively. In
Definition \ref{1.2.18'} we introduce the remaining five $\cat$-valued nerves of bicategories.
These associate to every bicategory $\c$ simplicial objects in $\cat$:  Those called its {\em
unitary categorical geometric nerves}, which are denoted by $\cgner\hspace{-4pt}^{^\mathrm{u}}\c$
and $\overline{\nabla}\hspace{-2.5pt}_{^{_\mathrm{u}}}\c$,  its {\em categorical geometric nerves},
denoted by $\cgner\c$ and $\overline{\nabla}\c$, and its {\em Segal nerve}, denoted by
$\underline{\mathrm{S}}\c$, respectively. As we will prove,  each of these five simplicial
categories can be thought of as a `rectification' of the pseudo-simplicial nerve $\cner\c$
(\ref{ps1.1}) of the bicategory, since they model the same homotopy type as $\class\c$.
Furthermore, all these simplicial categories are `special' in the sense \cite{segal74, thomason}
that the Segal projection maps on them induce homotopy equivalences on classifying spaces;
moreover, the latter, $\underline{\mathrm{S}}\c$, is a weak 2-category in the sense of Tamsamani
\cite{tam} and Simpson \cite{sim}, that is, the Segal projection maps on it are surjective
equivalences of categories \cite{lack-paoli} (or see the proof of Theorem \ref{psin} below).
Unitary categorical geometric nerves for 2-categories were treated in \cite{b-c} and for monoidal
categories in \cite{b-c2}.  Segal nerves of bicategories, also called 2-nerves in
\cite{lack-paoli},  enjoy many interesting properties whose study the cited paper by Lack and Paoli
is mainly dedicated to.

Considering a lax functor $F:[p]\rightsquigarrow \c$ as a family of `paths'  $F_{i,j}:Fj\to Fi$
between objects of the bicategory $\c$, and a lax transformation $\alpha:F\Rightarrow F'$ as a
`free homotopy'  between them, it is natural to consider {\em lax transformations relative to
objects} (i.e., to end points) as those $\alpha$ that are constantly identities on objects, that
is, with $\alpha i=1_{F_i}=1_{F'_i}$, for all $0\leq i\leq p$.  This kind of lax transformations
has been considered in \cite{b-c}, and also in \cite{lack, lack-paoli, g-g} under the name of {\em
icons} (short for `identity component oplax natural transformations').

\begin{remark} \label{monf}{\em
If $\m$ and $\m'$ are monoidal categories, then a monoidal functor $F:\m\to \m'$ amounts precisely
to a homomorphism between the corresponding delooping bicategories
$\Omega^{^{-1}}\hspace{-5pt}F:\Omega^{^{-1}}\hspace{-5pt}\m\rightsquigarrow
\Omega^{^{-1}}\hspace{-5pt}\m'$ (see Remark \ref{mono}). If $F,F':\m\to \m'$ are monoidal functors,
then a monoidal transformation $\alpha:F\Rightarrow F'$ is the same thing as a  relative to objects
lax transformation
$\Omega^{^{-1}}\hspace{-5pt}\alpha:\Omega^{^{-1}}\hspace{-2pt}F\Rightarrow\Omega^{^{-1}}\hspace{-2pt}F'$.
However, note that an arbitrary lax transformation $\beta:\Omega^{^{-1}}\hspace{-2pt}F\Rightarrow
\Omega^{^{-1}}\hspace{-2pt}F'$ consists of an object $x_0$ of $\m'$ (=$\beta *$, the component at
the unique object of $\Omega^{^{-1}}\hspace{-5pt}\m$) and morphisms $\widehat{\beta}_x: x_0\otimes
Fx\to F'\!x\otimes x_0$, one for each object $x$ of $\m$, satisfying the usual two conditions (see
Section 2).}\qed
\end{remark}

\begin{definition}\label{1.2.18'} The {\em categorical geometric nerves} $\cgner\c$ and $\overline{\nabla}\c$ of a bicategory
${\c}$ are the simplicial categories
$$ \begin{array}{ll}\cgner\c:&\Delta^{\!^{\mathrm{op}}}\ \to \ \cat,\\[6pt]
&[p]\mapsto \clfunc([p],\c)\end{array}$$ whose category of $p$-simplices is the category of  lax
functors $F:[p]\rightsquigarrow \c$, with relative to objects lax transformations between them as
arrows and
$$\begin{array}{ll}\overline{\nabla}\c&\Delta^{\!^{\mathrm{op}}}\ \to \ \cat,\\[6pt]
&[p]\mapsto \colfunc([p],\c)\end{array}$$ whose category of $p$-simplices is the category of  oplax
functors $F:[p]\rightsquigarrow \c$, with relative to objects lax transformations between them as
arrows, respectively.

The {\em unitary categorical geometric nerves}
 $\cgner\hspace{-4pt}^{^\mathrm{u}}\c \subseteq \cgner\c$ and
 $\overline{\nabla}\hspace{-2.5pt}_{^{_\mathrm{u}}}\c\subseteq \overline{\nabla}\c$ are the
 respective simplicial subcategories of normal lax functors and normal oplax functors $F:[p]\rightsquigarrow\c$, that is,
$$ \begin{array}{ll}\cgner\hspace{-4pt}^{^\mathrm{u}}\c:&\Delta^{\!^{\mathrm{op}}}\ \to \ \cat,\\[6pt]
&[p]\mapsto \cnlfunc([p],\c)\end{array}$$ and

$$ \begin{array}{ll}\overline{\nabla}\hspace{-2.5pt}_{^{_\mathrm{u}}}\c:&\Delta^{\!^{\mathrm{op}}}\ \to \ \cat.\\[6pt]
&[p]\mapsto \conlfunc([p],\c)\end{array}$$

The {\em Segal nerve} $\underline{\mathrm{S}}\c\subseteq \cgner\hspace{-4pt}^{^\mathrm{u}}\c$ is
the simplicial subcategory of normal homomorphisms $F:[p]\rightsquigarrow\c$, that is,

$$ \begin{array}{ll} \underline{\mathrm{S}}\c:&\Delta^{\!^{\mathrm{op}}}\ \to \ \cat.\\[6pt]
&[p]\mapsto \nhom([p],\c)\end{array}$$
\end{definition}

Note that, for $F,F':[p]\rightsquigarrow\c$ two geometric $p$-simplices, the existence of morphisms
$\widehat{\alpha}:F \Rightarrow F'$ in the categorical geometric nerve requires that $Fi=F'i$ for
all $0\leq i\leq p$, and, in such a case, an $\widehat{\alpha}:F\Rightarrow F'$ consists of a
family of deformations   in $\c$ \begin{equation}\label{1.2.21'} \xymatrix@C=9pt@R=9pt{ Fj
\ar@/^1pc/[rr]^{\textstyle F_{i,j}} \ar@/_1pc/[rr]_{\textstyle F'_{i,j}} & \scriptstyle{\textstyle
\Downarrow \widehat{\alpha}_{i,j}} & Fi,}
\end{equation}
for $0\leq i \leq j\leq p$, such that
\begin{equation}\label{1.2.18}
\xymatrix@C=9pt@R=9pt{ && Fj \ar@/_1pc/[lldd]_{\textstyle F_{i,j}} \ar@/^1pc/[lldd]^{\textstyle
F'_{i,j}} &&& &&& Fj
\ar@/_1pc/[lldd]_{\textstyle F_{i,j}} && \\
&\Downarrow \widehat{\alpha}&& \Downarrow\widehat{\alpha}&&&&& \ar@{}[d]^>>>>>>>{\textstyle\Downarrow\! \widehat{F}}&&\\
Fi &&\ar@{}[d]^{\textstyle \Downarrow\widehat{F'}}&& Fk \ar@/_1pc/[lluu]_{\textstyle F_{j,k}}
\ar@/^1pc/[lluu]^{\textstyle F'_{j,k}} \ar@/^3pc/[llll]^{\textstyle F'_{i,k}} & = & Fi &&&& Fk
\ar[llll]|<<<<<<<<<<{\textstyle F_{i,k}} \ar@/^3pc/[llll]^{\textstyle F'_{i,k}}\ar@/_1pc/[lluu]_{\textstyle F_{j,k}} \\
&&  &&&&&&\ \  \Downarrow\!\widehat{\alpha}}
\end{equation}
for  $0\leq i\leq j\leq k\leq p$, and, for $0\leq i\leq p$,
$$
\xymatrix{ &&&&&&\\Fi \ar@/^2.5pc/[rr]^{\textstyle 1_{Fi}} \ar[rr]|{\textstyle
F_{i,i}}\ar@/_2.5pc/[rr]_{\textstyle F'_{i,i}}\ar@{}[rru]|{\textstyle \Downarrow \widehat{F}}
\ar@{}[rrd]|{\textstyle \Downarrow \widehat{\alpha}} &  & Fi&=&Fi \ar@/^2.5pc/[rr]^{\textstyle
1_{Fi}} \ar@/_2.5pc/[rr]_{\textstyle F'_{i,i}}\ar@{}[rr]|{\textstyle \Downarrow \widehat{F'}}
 &  & Fi.\\
&&&&&&&}
$$

We now  complete the description of diagram (\ref{dia}): \vspace{0.2cm}
\begin{definition}\label{defdia}
\begin{enumerate}
\item[(i)] The simplicial maps $\gner\hspace{-4pt}^{^\mathrm{u}}\c \hookrightarrow \gner\c$ and $\nabla\hspace{-2.5pt}_{^{_\mathrm{u}}}\c\hookrightarrow \nabla\c$ are inclusions.

    \vspace{0.2cm}
\item[(ii)] The simplicial functors $\underline{\mathrm{S}}\c\hookrightarrow \cgner\hspace{-4pt}^{^\mathrm{u}}\c \hookrightarrow \cgner\c$ and $\overline{\nabla}\hspace{-2.5pt}_{^{_\mathrm{u}}}\c\hookrightarrow \overline{\nabla}\c$ are all inclusions.

\vspace{0.2cm}

\item[(iii)]The simplicial functor $\underline{\mathrm{S}}\c\hookrightarrow \overline{\nabla}\hspace{-2.5pt}_{^{_\mathrm{u}}}\c$  is the full embedding defined by
$$\left\{\begin{array}{l} (F,\widehat{F})\ \mapsto\ (F,\widehat{F}^{-1}), \\[4pt]
(F,\widehat{F})\overset{\widehat{\alpha}}\Longrightarrow (F',\widehat{F'})\ \mapsto \
(F,\widehat{F}^{-1})\overset{\widehat{\alpha}}\Longrightarrow
(F',\widehat{F'}^{-1}).\end{array}\right.
$$

    \vspace{0.2cm}

\item[(iv)]Each geometric nerve is regarded as a discrete simplicial category (i.e., with only identities) and the simplicial functors
$$\gner\c \hookrightarrow \cgner\c,\  \gner\hspace{-4pt}^{^\mathrm{u}}\c\hookrightarrow \cgner\hspace{-4pt}^{^\mathrm{u}}\c,\
\nabla\c\hookrightarrow\overline{\nabla}\c,\
\nabla\hspace{-2.5pt}_{^{_\mathrm{u}}}\c\hookrightarrow
\overline{\nabla}\hspace{-2.5pt}_{^{_\mathrm{u}}}\c   $$ are all inclusions, as each geometric
nerve is the simplicial set of objects of the corresponding categorical geometric nerve.

    \vspace{0.2cm}

\item[(v)] The pseudo-simplicial functor $J:\cner\c\rightsquigarrow \underline{\mathrm{S}}\c$  is defined as follows:

\end{enumerate}

\end{definition}
First, note that we have the equalities
$$
\begin{array}{cllll}
 \underline{\mathrm{S}}\c_0 & = & \mbox{Ob}\c= \cner\c_0,& &  \\
  &&&&\\
  \underline{\mathrm{S}}\c_1 & = & \bigsqcup\limits_{x_0,x_1}{\c}(x_1,x_0) & = & \cner\c_1 \\
\end{array}
$$
by identifying a normal homomorphism $F:[0]\rightsquigarrow\c$ with the object $F0$ and a normal
homomorphism $F:[1]\rightsquigarrow\c$ with the arrow $F_{0,1}:F1\to F0$ (see Lemma \ref{gn.1}). We
then take  $J_n:\cner\c_n\to \underline{\mathrm{S}}\c_n$  to be the identity functor for $n=0,1$.

For each $n\geq 2$, let
$$
J_n: \cner\c_n \hookrightarrow \underline{\mathrm{S}}\c_n
$$
be the functor taking an object of $\cner\c_n= \bigsqcup\limits _{(x_0,\ldots,x_n)\in
\mbox{\scriptsize Ob}\c^{n+1}}\hspace{-0.6cm} \c(x_1,x_0)\times\cdots\times\c(x_n,x_{n-1})$, say
the one given by the string
$$
\xymatrix@C=25pt{F0&\ar[l]_-{\textstyle F_{0,1}}F1&\ar[l]_{\textstyle F_{1,2}}\cdots&\ar[l]_{
\textstyle F_{n-1,n}}Fn}
$$
to the normal geometric $n$-simplex of $\c$ \begin{equation}\label{fe} F^{^{\mathrm
e}}=(F^{^{\mathrm e}},\widehat{F}^{^{\mathrm e}}):[n]\rightsquigarrow \c ,\end{equation} which
places the objects $F^{^{\mathrm e}}\!i=Fi$ at the vertices and the morphism $F^{^{\mathrm
e}}_{i,i+1}=F_{i,i+1}:F_{i+1}\rightarrow F_i$ at the edges $i+1\to i$. For $0\leq i<j<n$, the
arrows  $F^{^{\mathrm e}}_{i,j+1}:F_{j+1}\to F_i$ are then inductively given by
$$ F^{^{\mathrm e}}_{i,j+1}=F^{^{\mathrm e}}_{i,j}\circ F_{j,j+1},$$
so that the triangles
$$
\xymatrix@C=25pt@R=25pt{ &\hspace{6pt} Fj \hspace{4pt}\ar[dl]_{\textstyle F^{^{\mathrm e}}_{i,j}}            \\
 Fi  & &  \ar[ll]^{\textstyle F^{^{\mathrm e}}_{i,j+1}}    F(j+1)\ar[ul]_{\textstyle F^{^{\mathrm e}}_{j,j+1}}, }
$$
commute and the deformations $\widehat{F}^{^{\mathrm e}}_{i,j,j+1}:{F}^{^{\mathrm e}}_{i,j}\circ
{F}^{^{\mathrm e}}_{j,j+1}\Rightarrow {F}^{^{\mathrm e}}_{i,j+1}$ are all taken to be identities.
For $0\leq i<j<k< n$, the 2-cells $\widehat{F}^{^{\mathrm e}}_{i,j,k+1}:{F}^{^{\mathrm
e}}_{i,j}\circ {F}^{^{\mathrm e}}_{j,k+1}\Rightarrow {F}^{^{\mathrm e}}_{i,k+1}$ are canonically
determined by the associativity constraints of $\c$ through the commutative diagrams
$$
\xymatrix@C=45pt{ {F}^{^{\mathrm e}}_{i,j}\circ {F}^{^{\mathrm e}}_{j,k+1} \ar@2{->}[r]^{\textstyle \widehat{F}^{^{\mathrm e}}_{i,j,k+1}} & {F}^{^{\mathrm e}}_{i,k+1}\ar@{=}[d]\\
({F}^{^{\mathrm e}}_{i,j}\circ F^{^{\mathrm e}}_{j,k})\circ F^{^{\mathrm e}}_{k,k+1}
\ar@2{->}[u]^{\textstyle \boldsymbol{a}}\ar@2{->}[r]^-{\textstyle \widehat{F}^{^{\mathrm
e}}_{i,j,k}\circ 1} & {F}^{^{\mathrm e}}_{i,k}\circ F^{^{\mathrm e}}_{k,k+1}. }
$$

Thus,  for example:

\noindent$\xymatrix@R=5pt@C=15pt{&&F1\ar[ddl]_{F_{0,1}}\ar@{}[dd]|>>>>>>{\textstyle \Downarrow 1}&\\  J_2(F0\overset{F_{0,1}}\longleftarrow F1\overset{F_{1,2}}\longleftarrow F2)=&&& \\
&F0&&\ar[uul]_{F_{1,2}} F2\,,\ar[ll]^{F_{0,1}\circ F_{1,2}}}$

\noindent$ \xymatrix { &
 & F3 \ar[dl]_{(F_{0,1}\circ F_{1,2})\circ F_{2,3}}  \ar[dd]|<<<<<<{F_{1,2}\circ F_{2,3} }\ar[dr]^{ F_{2,3}}& \\
 J_3(F0\overset{F_{0,1}}\longleftarrow F1\overset{F_{1,2}}\longleftarrow F2\overset{F_{2,3}}\longleftarrow F3)=&
F0   & & \text{-}F2 \ar[dl]^{ F_{1,2}} \ar[ll]|<<<<<<{\; F_{0,1}\circ F_{1,2}\;}\\
 && F1 \ar[ul]^{
F_{0,1}}} $

\noindent with the deformations
$$\begin{array}{l}
  \widehat{F}^{^{\mathrm e}}_{0,1,2}=1:F_{0,1}\circ F_{1,2}\Rightarrow F_{0,1}\circ F_{1,2}\\
\widehat{F}^{^{\mathrm e}}_{0,1,3}= \boldsymbol{a}^{-1}:F_{0,1}\circ (F_{1,2}\circ F_{2,3})\Rightarrow (F_{0,1}\circ F_{1,2})\circ F_{2,3}\\
\widehat{F}^{^{\mathrm e}}_{0,2,3}=1:(F_{0,1}\circ F_{1,2})\circ F_{2,3}\Rightarrow (F_{0,1}\circ F_{1,2})\circ F_{2,3}  \\
\widehat{F}^{^{\mathrm e}}_{1,2,3}=1:F_{1,2}\circ F_{2,3}\Rightarrow F_{1,2}\circ F_{2,3},
\end{array}
$$
and so on.

It is a consequence of the coherence theorem for bicategories that each
$$J_n(F_0\overset{F_{0,1}}\longleftarrow F_1 \overset{F_{1,2}}\longleftarrow \cdots
\overset{F_{n-1,n}} \longleftarrow F_n)=(F^{^{\mathrm e}},\widehat{F}^{^{\mathrm
e}}):[n]\rightsquigarrow \c$$ is in fact a normal geometric $n$-simplex of $\c$ (recall Lemma
\ref{gn.1}). Indeed, since every deformation $\widehat{F}^{^{\mathrm e}}_{i,j,k}:{F}^{^{\mathrm
e}}_{i,j}\circ {F}^{^{\mathrm e}}_{j,k}\Rightarrow {F}_{i,k}$ is invertible, it is actually  a
normal homomorphism  $F^{^{\mathrm e}}:[n]\rightsquigarrow\c$.

Further, the functor $J_n$ on an arrow
$$
\xymatrix @C=1pc {F0  & {\Downarrow\, \alpha_{0,1}} & F1 \ar@/^1pc/[ll]^{\textstyle G_{0,1}}
\ar@/_1pc/[ll]_-{\textstyle F_{0,1}} & {\Downarrow\, \alpha_{1,2}} & F2  \ar@/^1pc/[ll]^{\textstyle
G_{1,2}} \ar@/_1pc/[ll]_-{\textstyle F_{1,2}}&\hspace{-15pt}\cdots&\hspace{-10pt}
 F(n-1)
  & {\Downarrow\, \alpha_{n-1,n}} & Fn
\ar@/^1pc/[ll]^{\textstyle G_{n-1,n}} \ar@/_1pc/[ll]_-{\textstyle F_{n-1,n}}}
$$
in $\cner\c_n$ is the arrow in $\underline{\mathrm{S}}\c_n$ (i.e., the relative to objects lax
transformation) ${\widehat{\alpha}^{_{\mathrm e}}:F^{^{\mathrm e}}\Rightarrow G^{^{\mathrm e}}}$,
as in (\ref{1.2.21'}),  inductively given by the deformations
$$\widehat{\alpha}^{_{\mathrm e}}_{i,j+1}=\left\{%
\begin{array}{lll}
    \alpha_{i,i+1} & \text{if} & j=i \\
    &&\\
\widehat{\alpha}^{_{\mathrm e}}_{i,j}\circ \alpha_{j,j+1} & \text{if} & j>i. \\
\end{array}%
\right.$$

From its construction, it is clear that each $J_n:\cner\c_n \rightarrow \underline{\mathrm{S}}\c_n$
is a faithful and injective on objects functor. That $J_n$ is full follows from the equalities
(\ref{1.2.18}) since every geometric simplex $F:[n]\rightsquigarrow \c$ in
$\underline{\mathrm{S}}\c$ is a homomorphism, that is, the  deformations $\widehat{F}_{i,j,k}$ are
all invertible.

The family of functors $J_n:\cner\c_n\to \underline{\mathrm{S}}\c_n$, together with the natural
isomorphisms $J_{n-1}d_i\cong d_iJ_n$ and $J_ns_i\cong s_i J_{n-1}$, canonically induced by the
associativity and unit constraints $\boldsymbol{a}$, $\boldsymbol{l}$, and $\boldsymbol{r}$ of the
bicategory $\c$, becomes a morphism of supercoherent structures \cite{jardine}, and this defines
our pseudo-simplicial full embedding $J:\cner\c\rightsquigarrow \underline{\mathrm{S}}\c$.
Actually, $J$ is a pseudo equivalence, as we will show later in the proof of Theorem \ref{psin}.

To end this section, we stress that the construction  $\cner\c$ is functorial on homomorphisms
between bicategories $\c\rightsquigarrow\c'$, whereas  $\underline{\mathrm{S}}\c$ is on normal
homomorphisms, $\gner\hspace{-4pt}^{^\mathrm{u}}\c$ and  $\cgner\hspace{-4pt}^{^\mathrm{u}}\c$ are
on normal lax functors, $\gner\c$ and $\cgner\c$ are on lax functors,
$\nabla\hspace{-2.5pt}_{^{_\mathrm{u}}}\c$ and
$\overline{\nabla}\hspace{-2.5pt}_{^{_\mathrm{u}}}\c$ are on normal oplax functors, and  $\nabla\c$
and $\overline{\nabla}\c$ are on all oplax functors.

\section{The homotopy invariance theorem} In this section, our goal is to prove the following:
\begin{theorem}\label{thit}
 For any bicategory $\c$, all the maps in the diagram
$$
\xymatrix{&&\ \  \class\c \ar@<1.5pt>[d]\ar@<1.5pt>[dl]\ \ \ar@<-1.5pt>[dr]&&\\ \class\cgner\c\ &
\ar@<1.5pt>[l]\ \class\cgner\hspace{-3pt}^{^\mathrm{u}}\c\ &  \ar@<1.5pt>[l]\
\class\underline{\mathrm{S}}\c \ \ar@<-1.5pt>[r] &\
\class\overline{\nabla}\hspace{-2.5pt}_{^{_\mathrm{u}}}\c\  \ar@<-1.5pt>[r]& \
\class\overline{\nabla}\c
\\ \ar@<1.5pt>[u]\class\gner\c\
&\ar@<1.5pt>[l]\ \class\gner\hspace{-3pt}^{^\mathrm{u}}\c \ \ar@<1.5pt>[u] &&\
\class\nabla\hspace{-2.5pt}_{^{_\mathrm{u}}}\c \ \ar@<1.5pt>[u] \ar@<-1.5pt>[r] & \ \class\nabla\c,
\ \ar@<1.5pt>[u]}
$$
induced by (\ref{dia}) on classifying spaces, are homotopy equivalences.
\end{theorem}
By symmetry, we only will attack the problem of proving what concerns the maps in the middle and on
the left of the diagram. The other homotopy equivalences are proven in a parallel way. The proof
can be  naturally divided into four theorems, given  below.

\begin{theorem}For any bicategory $\c$, both simplicial inclusion functors ${\gner\c\hookrightarrow \cgner\c}$  and $\gner\hspace{-3pt}^{^\mathrm{u}}\c\hookrightarrow\cgner\hspace{-3pt}^{^\mathrm{u}}\c$ induce
homotopy equivalences on the corresponding classifying spaces, thus $ \class\gner\c \simeq
\class\cgner\c$ and $\class\gner\hspace{-3pt}^{^\mathrm{u}}\c\simeq
\class\cgner\hspace{-3pt}^{^\mathrm{u}}\c$.
\end{theorem}
\begin{proof}

The following argument  is for ${\gner\c\hookrightarrow \cgner\c}$, and we leave it to the reader
to check that all the following constructions below restrict to the subinclusion
$\gner\hspace{-3pt}^{^\mathrm{u}}\c\hookrightarrow\cgner\hspace{-3pt}^{^\mathrm{u}}\c$, there by
obtaining a corresponding proof for this latter.

Let  $\ner\cgner\c$ be the bisimplicial set obtained from the simplicial category $\cgner\c$ by
composing with the nerve of categories functor (\ref{p.1.1}). Then, a $(p,q)$-simplex of
$\ner\cgner\c$  is  a string
$$
\xi=(F^0\overset{\ \textstyle \widehat{\alpha}^1}\Longleftarrow F^1 \overset{\ \textstyle
\widehat{\alpha}^2}\Longleftarrow \cdots \overset{\ \textstyle \widehat{\alpha}^p}\Longleftarrow
F^p)
$$
of $p$ composible arrows in the category  $\cgner\c_q=\clfunc([q],\c)$, of lax functors
$[q]\rightsquigarrow\c$. The vertical face and degeneracy operators in $\ner\cgner\c$
$$
\xymatrix@C=9pt@R=9pt{ \ner\cgner\c_{p,q+1} && \ner\cgner\c_{p,q} \ar[ll]_-{\textstyle s_m^v}
\ar[rr]^-{\textstyle d_m^v} && \ner\cgner\c_{p,q-1}, && 0\leq m\leq q,\\ }
$$
 are induced by those of $\cgner\c$, that is,
$$\xymatrix{d_m^v\xi=(d_mF^0&\ar@{=>}[l]_-{\textstyle d_m\widehat{\alpha}^1}d_mF^1& \ar@{=>}[l]_-{ \textstyle d_m\widehat{\alpha}^2} \cdots& \ar@{=>}[l]_-{ \textstyle d_m\widehat{\alpha}^p} d_mF^p),}$$
and similarly  $s_m^v\xi$;  the  horizontal ones
$$
\xymatrix@C=9pt@R=9pt{ \ner\cgner\c_{p+1,q} && \ner\cgner\c_{p,q} \ar[ll]_-{\textstyle s_m^h}
\ar[rr]^-{\textstyle d_m^h} && \ner\cgner\c_{p-1,q}, && 0\leq m\leq p,\\ }
$$
are those of  $\ner(\cgner\c_q)$, which are given by the formulas (\ref{p.1.2}).

Since $\gner\c$ is a discrete simplicial category,  $\ner\gner\c$ is a bisimplicial set that is
constant in the horizontal direction. The induced bisimplicial inclusion
$\ner\gner\c\to\ner\cgner\c$ is then, for each $p\geq 0$, the composite simplicial map
\begin{equation}\label{1.2.28'}\begin{array}{l} \gner\c=\ner\cgner\c_{0,*}\overset{s_0^h}\hookrightarrow \ner\cgner\c_{1,*}\overset{s_0^h}\hookrightarrow  \cdots \overset{s_0^h}\hookrightarrow \ner\cgner\c_{p,*}.\\[5pt]
\xymatrix{F\ar@{|->}[r]& (F\overset{\textstyle \widehat{1}}\Longleftarrow \underset{~}F
\overset{\textstyle \widehat{1}}\Longleftarrow \overset{(p}\cdots \overset{\ \textstyle
\widehat{1}}\Longleftarrow F)}\end{array}\end{equation}

Taking into account the homeomorphisms (\ref{p.1.9}) to prove that $\class\gner\c\to\class\cgner\c$
is a homotopy equivalence, we are going to  prove that the induced simplicial map on diagonals $
\gner\c\to \diag \ner\cgner\c$ is a  weak equivalence. To do so, as every pointwise weak homotopy
equivalence bisimplicial map is a diagonal  weak homotopy equivalence \cite[IV, Proposition
1.7]{g-j}, it suffices to  prove that  each of these simplicial maps (\ref{1.2.28'}) is a weak
homotopy equivalence. Actually, we will prove more: {\em Every simplicial map $$\begin{array}{ll}
s_0^h:&\ner\cgner\c_{p,*}\hookrightarrow \ner\cgner\c_{p+1,*},\\[4pt]
&(F^0\overset{\ \textstyle \widehat{\alpha}^1}\Longleftarrow F^1 \overset{\ \textstyle
\widehat{\alpha}^2}\Longleftarrow \cdots \overset{\ \textstyle \widehat{\alpha}^p}\Longleftarrow
F^p)\mapsto (F^0\overset{\textstyle \widehat{1}}\Longleftarrow F^0\overset{\ \textstyle
\widehat{\alpha}^1}\Longleftarrow F^1 \overset{\ \textstyle \widehat{\alpha}^2}\Longleftarrow
\cdots \overset{\ \textstyle \widehat{\alpha}^p}\Longleftarrow F^p)\end{array}
$$
 embeds the simplicial set $\ner\cgner\c_{p,*}$ into $\ner\cgner\c_{p+1,*}$ as a simplicial deformation retract.}

To do so, since $d_0^hs_0^h=1$, it is enough to exhibit a simplicial homotopy
$$
H^p:1\Rightarrow s_0^hd_0^h:\ner\cgner\c_{p,*}\to\ner\cgner\c_{p,*},
$$
for each $p\geq 1$.

We first consider the case $p=1$. Then, to define $H^1$, we begin by constructing a simplicial
homotopy
\begin{equation}\label{1.2.28}
T:d_1^h\Rightarrow d_0^h: \ner\cgner\c_{1,*}\rightarrow \ner\cgner\c_{0,*}=\gner\c,
\end{equation}
$$ \xymatrix{\cdots & \ner\cgner\c_{1,q+1} \ar@<1.5ex>[rr]^{d_0^v} \ar@<-1.5ex>[rr]_{d_{q+1}^v}
\ar@<1ex>[dd]^{d_0^h} \ar@<-1ex>[dd]_{d_1^h} & {}^{^{\vdots}} &
\ner\cgner\c_{1,q}\ar@<1ex>[dd]^{d_0^h} \ar@<-1ex>[dd]_{d_1^h}
\ar@<1ex>[ddll]^{T_0}_{\ \ \dots}\ar@<-1ex>[ddll]_{T_q} & \cdots \\
&&&&\\
\cdots & \gner\c_{q+1} \ar@<1.5ex>[rr]^{d_0} \ar@<-1.5ex>[rr]_{d_{q+1}} &{}^{^{\vdots}}& \gner\c_q
&
\cdots\\
}
$$
by the maps $T_m:\ner\cgner\c_{1,q}\to \gner\c_{q+1}$, $0\leq m\leq q$, defined as follows: For any
arrow $\widehat{\alpha}:F^1\Rightarrow F^0$ in the category $\cgner\c_q=\clfunc([q],\c)$, and
writing $F^1j=F^0j$ by $Fj$  for $0\leq j\leq q$, then
$$T_m(\widehat{\alpha}):[q+1]\rightsquigarrow\c$$is  the geometric $(q+1)$-simplex of $\c$ that
places

\vspace{0.2cm}
\noindent - at each vertex $i$, $0\leq i\leq q+1$, the object $\left\{%
\begin{array}{lll}
    Fi & \text{if} & i\leq m, \\
    F(i-1) & \text{if} & i>m; \\
\end{array}%
\right. $

\noindent - at each edge $j\to i$, $0\leq i \leq j\leq q+1$, the arrow

$$\left\{
\begin{array}{lll}
    F_{i,j}^1:Fj\to Fi & \text{if} & j\leq m, \\
    &&\\
    F_{i,j-1}^0:F(j-1)\to Fi & \text{if} & i\leq m <j, \\
    &&\\
    F_{i-1,j-1}^0:F(j-1)\to F(i-1) & \text{if} &  m<i; \\
\end{array}\right.$$

\noindent \hspace{-0.5cm}$\xymatrix@C=8pt@R=0pt{ & &j\ar[ddl] & &\\\text{- at each triangle}& &&,&
0\leq i \leq j \leq k\leq q+1, \text{the 2-cell in }\c\\  &   i && k\ar[ll]\ar[uul]&}$

$$\left\{\begin{array}{l}
\xymatrix{\ar@{}[drr]|(.6)*+{\Downarrow \widehat{F}^1_{_{i,j,k}}}
                & Fj \ar[dl]_{\textstyle F^1_{i,j}}             \\
Fi  & &     Fk\,,   \ar[ul]_{\textstyle F^1_{j,k}} \ar[ll]^{\textstyle F^1_{i,k}}&\hspace{0.6cm}
\text{if} \hspace{0.6cm} k\leq m,&   }
\\
\xymatrix@C=50pt@R=50pt{\ar@{}[drr]|(0.7)*+{\Downarrow \widehat{F}^0_{_{i,j,k-1}}}
                & Fj \ar[dl]_{\textstyle F^1_{i,j}}^>>>>>>>>>>>>>{\textstyle\Downarrow \!\widehat{\alpha}_{_{i,j}}}   \ar@/^28pt/[dl]^<<<<<<<{\textstyle F_{i,j}^0}          \\
Fi  & &     F(k-1)\,,   \ar[ul]_{\textstyle F^0_{j,k-1}} \ar[ll]^{\textstyle
F^0_{i,k-1}}&\hspace{-1cm}   \text{ if }\hspace{0.6cm}  j\leq m< k, &  }
\\
\xymatrix{\ar@{}[drr]|(.6)*+{\Downarrow \widehat{F}^0_{_{i,j-1,k-1}}}
                & F(j-1) \ar[dl]_{\textstyle F^0_{i,j-1}}             \\
Fi  & &     F(k-1)\,,   \ar[ul]_{\textstyle F^0_{j-1,k-1}} \ar[ll]^{\textstyle
F^0_{i,k-1}}&\hspace{0.6cm}   \text{if} \hspace{0.6cm} i\leq m< j,&   }
\\

\xymatrix{\ar@{}[drr]|(.6)*+{\Downarrow \widehat{F}^0_{_{i-1,j-1,k-1}}}
                & F(j-1) \ar[dl]_{\textstyle F^0_{i-1,j-1}}             \\
F(i-1)  & &     F(k-1)\,,   \ar[ul]_{\textstyle F^0_{j-1,k-1}} \ar[ll]^{\textstyle
F^0_{i-1,k-1}}&\hspace{0.6cm}   \text{if} \hspace{0.6cm} m< i;&   }
\end{array}\right.
$$

\noindent and whose unit structure deformations (\ref{1.2.16'}) are
$$
\left\{\begin{array}{lll}\widehat{F}^1_i:1_{Fi}\Rightarrow F^1_{i,i},& \text{if }& i\leq m,\\[6pt]
\widehat{F}^0_{i-1}:1_{F(i-1)}\Rightarrow F^0_{i-1,i-1},& \text{if } &m< i.\end{array}\right.
$$

\vspace{0.2cm}

It is straightforward to verify that the above definitions do indeed give the  simplicial homotopy
(\ref{1.2.28}) as predicted. Further, note that we have the equalities $T_0=s_0d_1^h$ and
$T_ms_0^h=s_m$, for all $0\leq m \leq q$.

Then, the simplicial homotopy $H^1:1\Rightarrow  s_0^hd_0^h$
$$ \xymatrix{\cdots & \ner\cgner\c_{1,q+1} \ar@<1.5ex>[rr]^{d_0^v} \ar@<-1.5ex>[rr]_{d_{q+1}^v}
\ar@<1ex>[dd]^{s_0^hd_0^h} \ar@<-1ex>[dd]_{1} & {}^{^{\vdots}} &
\ner\cgner\c_{1,q}\ar@<1ex>[dd]^{s_0^hd_0^h} \ar@<-1ex>[dd]_{1}
\ar@<1ex>[ddll]^{H_0^1}_{\ \ \dots}\ar@<-1ex>[ddll]_{H_q^1} & \cdots \\
&&&&\\
\cdots & \ner\cgner\c_{1,q+1} \ar@<1.5ex>[rr]^{d_0^v} \ar@<-1.5ex>[rr]_{d_{q+1}^v}
&{}^{^{\vdots}}&\ner \cgner\c_{1,q} &
\cdots\\
}
$$
is defined by the maps $H^1_m$, $0\leq m \leq q$, as above, which apply an arrow
$\widehat{\alpha}:F^1\Rightarrow F^0$ in $\ner\cgner\c_q$ to the arrow
\begin{equation}\label{hhat} \xymatrix{s_mF^1\ar@{=>}[r]^{H_m^1(\widehat{\alpha})}& T_m(\widehat{\alpha})}\end{equation}
in $\ner\cgner\c_{q+1}$, given by the family of deformations of $\c$
$$
H_m^1(\widehat{\alpha})_{i,j}=\left\{%
\begin{array}{lll}
    1:F_{i,j}^1\Rightarrow F_{i,j}^1& \text{if} & j\leq m, \\[4pt]
    \alpha_{i,j-1}:F_{i,j-1}^1\Rightarrow F_{i,j-1}^0 & \text{if} & i\leq m<j, \\[4pt]
     \alpha_{i-1,j-1}:F_{i-1,j-1}^1\Rightarrow F_{i-1,j-1}^0 & \text{if} & i<m, \\
\end{array}%
\right.
$$
for $0\leq i\leq j\leq q+1$.

One easily observes that $H_0^1(\widehat{\alpha})=s_0^v(\widehat{\alpha})$ and that
${d_{q+1}^vH_q^1(\widehat{\alpha})=1_{F^1}}$, whence $d_0^vH_0^1=1$ and
$d_{q+^1}^vH_q^1=s_0^hd_0^h$.  Checking in full the remaining homotopy simplicial identities needed
for $H^1:1\Rightarrow s_0^hd_0^h$ to be a simplicial homotopy is again straightforward (though
quite tedious) and we leave it to the reader.

Finally, for an arbitrary $p\geq 2$, the simplicial homotopy $H^p:1\Rightarrow s_0^hd_0^h$
$$ \xymatrix{\cdots & \ner\cgner\c_{p,q+1} \ar@<1.5ex>[rr]^{d_0^v} \ar@<-1.5ex>[rr]_{d_{q+1}^v}
\ar@<1ex>[dd]^{s_0^hd_0^h} \ar@<-1ex>[dd]_{1} & {}^{^{\vdots}} &
\ner\cgner\c_{p,q}\ar@<1ex>[dd]^{s_0^hd_0^h} \ar@<-1ex>[dd]_{1}
\ar@<1ex>[ddll]^{H_0^p}_{\ \ \dots}\ar@<-1ex>[ddll]_{H_q^p} & \cdots \\
&&&&\\
\cdots & \ner\cgner\c_{p,q+1} \ar@<1.5ex>[rr]^{d_0^v} \ar@<-1.5ex>[rr]_{d_{q+1}^v}
&{}^{^{\vdots}}&\ner \cgner\c_{p,q} &
\cdots\\
}
$$
is given by the maps $H_m^p$, $0\leq m\leq q$, which take a $(p,q)$-simplex of $\ner\cgner\c$, say
$$\xi=(F^0\overset{\ \textstyle \widehat{\alpha}^1}\Longleftarrow F^1 \overset{\ \textstyle \widehat{\alpha}^2}\Longleftarrow \cdots \overset{\ \textstyle \widehat{\alpha}^p}\Longleftarrow F^p),$$
to the $(p,q+1)$-simplex
$$
H^p_m(\xi)=(\xymatrix@C=30pt{T_m(\widehat{\alpha}^1)&\ar@{=>}[l]_-{\textstyle
H^1_m(\widehat{\alpha}^1)} s_mF^1&\ar@{=>}[l]_-{\textstyle
s_m\widehat{\alpha}^2}s_mF^2\cdots&\ar@{=>}[l]_-{ \textstyle s_m\widehat{\alpha}^p}s_mF^p}),$$
where $H_m^1(\widehat{\alpha}^1)$ is as in (\ref{hhat}). The simplicial identities that make $H^p$
a simplicial homotopy are at this stage much easier to verify. This completes the proof.
\end{proof}

The proof of the following theorem involves the `normalization' process.

\begin{theorem}\label{1.2.24} For any bicategory $\c$, both the simplicial inclusion functor $\cgner\hspace{-3pt}^{^\mathrm{u}}\c\hookrightarrow \cgner\c$  and the simplicial inclusion map $\gner\hspace{-3pt}^{^\mathrm{u}}\c\hookrightarrow\gner\c$ induce
homotopy equivalences on the corresponding classifying spaces, and therefore
 $\class\cgner\hspace{-3pt}^{^\mathrm{u}}\c \simeq \class\cgner\c$ and
$\class\gner\hspace{-3pt}^{^\mathrm{u}}\c\simeq \class\gner\c$.
\end{theorem}
\begin{proof} Since we have the commutative square
$$
\xymatrix{\class\cgner\hspace{-3pt}^{^\mathrm{u}}\c\ar[r]&\class\cgner\c\\
\class\gner\hspace{-3pt}^{^\mathrm{u}}\c\ar[u]^-{\wr}\ar[r]&\class\gner\c,\ar[u]_-{\wr}}
$$
in which, by Theorem \ref{1.2.24} above,  both vertical maps are homotopy equivalences, it is
sufficient to prove that the simplicial functor
$\cgner\hspace{-3pt}^{^\mathrm{u}}\c\hookrightarrow \cgner\c$ induces a homotopy equivalence on
classifying spaces.

We will prove that, for each $p\geq 0$, the inclusion functor
$$
\cgner\hspace{-3pt}^{^\mathrm{u}}\c_p=\cnlfunc([p],\c)\ \hookrightarrow\ \cgner\c_p=\clfunc([p],\c)
$$
has a right adjoint. Then, every induced map $\class
\cgner\hspace{-3pt}^{^\mathrm{u}}\c_p\to\class\cgner\c_p$ is a homotopy equivalence and the result
follows from Theorem \ref{p.1.23}.

The right adjoint functor  $(-)\hspace{-2pt}^{^\mathrm{u}}:\cgner\c_p\to
\cgner\hspace{-3pt}^{^\mathrm{u}}\c_p$,  which should be called the {\em normalization functor},
works as follows: By Lemma \ref{gn.1}, every lax functor $F=(F,\widehat{F}):[p]\rightsquigarrow\c$
uniquely determines a normal lax functor
$F^{^\mathrm{u}}=(F^{^\mathrm{u}},\widehat{F}^{^\mathrm{u}}):[p]\rightsquigarrow \c$, such that
$$
\left\{\begin{array}{lll}F^{^\mathrm{u}}\!i=Fi& \text{for}& 0\leq i\leq p,\\[4pt]
F^{^\mathrm{u}}_{i,j}=F_{i,j}& \text{for}&0\leq i<j\leq p,\\[4pt]
\widehat{F}^{^\mathrm{u}}_{i,j,k}=\widehat{F}_{i,j,k}&\text{for}&0\leq i<j<k\leq
p.\end{array}\right.
$$
This normal lax functor $F^{^\mathrm{u}}$ is indeed the unique one such that the family of deformations in $\c$, $\widehat{\nu}_{i,j}=\left\{%
\begin{array}{lcl}
  1_{F_{i,j}}:F_{i,j}\Rightarrow F_{i,j} & \text{if} & i\neq j   \\
  \widehat{F}_i:1_{Fi}\Rightarrow F_{i,i} & \text{if} & i=j  \\
\end{array}%
 \right.$, $0\leq i \leq j\leq p$, becomes a morphism $\widehat{\nu}:F^{^\mathrm{u}}\Rightarrow F$ in $\cgner\hspace{-3pt}^{^\mathrm{u}}\c_p$.
Furthermore, any morphism $\alpha:F\Rightarrow G$ in $\cgner\c_p$ uniquely determines a morphism
$\alpha\hspace{-1pt}^{^{_\mathrm{u}}}:F^{^\mathrm{u}}\Rightarrow G^{^\mathrm{u}}$ in
$\cgner\hspace{-3pt}^{^\mathrm{u}}\c_p$ such that the square
$$\xymatrix@C=9pt@R=9pt{F^{^\mathrm{u}}\ar@2[rr]^{\textstyle\widehat{\nu}_{_F}} \ar@2[dd]_{\textstyle\alpha\hspace{-1pt}^{^{_\mathrm{u}}}} && F \ar@2[dd]^{\textstyle{\alpha}} \\
&&\\
G^{^\mathrm{u}}\ar@2[rr]^{\textstyle\widehat{\nu}_{_G}} && G\\ }$$ is commutative,
namely, $\alpha\hspace{-1pt}^{^{_\mathrm{u}}}_{i,j}=\left\{%
\begin{array}{lcl}
  \alpha_{i,j} & if & i\neq j   \\
  1_{1_{Fi}} & if & i=j  \\
\end{array}%
 \right.$. These mappings $F\mapsto F^{^\mathrm{u}}$, $\alpha \mapsto \alpha\hspace{-1pt}^{^{_\mathrm{u}}}$, describe
  the functor $(-)\hspace{-2pt}^{^\mathrm{u}}:\cgner\c_p\to \cgner\hspace{-3pt}^{^\mathrm{u}}\c_p$ that
   is right adjoint to the inclusion $\cgner\hspace{-3pt}^{^\mathrm{u}}\c_p\hookrightarrow \cgner\c_p$, with the identity
    and $\widehat{\nu}$ being the unit and the counit of the adjunction, respectively.
\end{proof}

In the proof of  Theorem \ref{psin} below, our discussion uses the so-called {\em Segal projection}
functors
 (see \cite[Definition 1.2]{segal74})
that, in our simplicial category $\cgner\hspace{-3pt}^{^\mathrm{u}}\c$, are defined by
$$
 P_n=\prod_{k=1}^n d_0\cdots d_{k-2} d_{k+1}\cdots d_n :\cgner\hspace{-3pt}^{^\mathrm{u}}\c_n\longrightarrow
\cner\c_n
$$
for $n\geq 2$. That is, each
$$
P_n:\cnlfunc([n],\c)\to\hspace{-0.3cm}\bigsqcup_{(x_0,\ldots,x_n)\in \mbox{\scriptsize
Ob}\c^{n+1}}\hspace{-0.5cm} \c(x_1,x_0)\times\c(x_2,x_1)\times\cdots\times\c(x_n,x_{n-1})
$$
maps a normal geometric $n$-simplex of $\c$, $F=(F,\widehat{F}):[n]\rightsquigarrow\c$ to the
string
$$
P_nF=(\xymatrix@C=25pt{F0&\ar[l]_-{\textstyle F_{0,1}}F1&\ar[l]_{\textstyle F_{1,2}}\cdots&\ar[l]_{
\textstyle F_{n-1,n}}Fn}),
$$
and a morphism  $\widehat{\alpha}:F\Rightarrow F'$ in $\cgner\hspace{-3pt}^{^\mathrm{u}}\c_n$ to
$$
\xymatrix @C=1pc {F0  & {\Downarrow\widehat{\alpha}_{0,1}} & F1 \ar@/^1pc/[ll]^{\textstyle
F'_{0,1}} \ar@/_1pc/[ll]_-{\textstyle F_{0,1}} & {\Downarrow \widehat{\alpha}_{1,2}} & F2
\ar@/^1pc/[ll]^{\textstyle F'_{1,2}} \ar@/_1pc/[ll]_-{\textstyle
F_{1,2}}&\hspace{-15pt}\cdots&\hspace{-10pt}
 F(n-1)
  & {\Downarrow\widehat{\alpha}_{n-1,n}} & Fn.
\ar@/^1pc/[ll]^{\textstyle F'_{n-1,n}} \ar@/_1pc/[ll]_-{\textstyle F_{n-1,n}}}
$$

\begin{theorem}\label{psin}
For any bicategory $\c$, the pseudo-simplicial embeddings $J:\cner\c \rightsquigarrow
\cgner\hspace{-3pt}^{^\mathrm{u}}\c$ and $J:\cner\c\rightsquigarrow \underline{\mathrm{ S}}\c$
induce both homotopy equivalences on the corresponding classifying spaces; thus, $ \class\c\simeq
\class\cgner\hspace{-3pt}^{^\mathrm{u}}\c $ and $\class\c\simeq \class\underline{\mathrm{ S}}\c$.
\end{theorem}
\begin{proof}
For $n=0,\ 1$, we have equalities
$$
\begin{array}{c}
\cgner\hspace{-3pt}^{^\mathrm{u}}\c_0=\underline{\mathrm{S}}\c_0  =  \mbox{Ob}\c= \cner\c_0  \\
  \\
\cgner\hspace{-3pt}^{^\mathrm{u}}\c_1=\underline{\mathrm{S}}\c_1  =  \bigsqcup\limits_{x_0,x_1}{\c}(x_1,x_0)  =  \cner\c_1 \\
\end{array}
$$
by identifying a normal lax functor $F:[0]\rightsquigarrow\c$ with the object $F0$, and a normal
lax functor $F:[1]\rightsquigarrow\c$ with the arrow $F_{0,1}:F1\to F0$. The functors $J_0$ and
$J_1$ are both identities.

For $n\geq 2$, we have the equality $P_nJ_n=1$. Moreover, there is a natural transformation
$\widehat{\varepsilon}:J_nP_n\Rightarrow 1$ defined as follows: If $F:[n]\rightsquigarrow \c$ is
any normal lax functor, that is, any object of the category
$\cgner\hspace{-3pt}^{^\mathrm{u}}\c_n$, and $J_nP_n(F)= F^{^{\mathrm e}}$ as in (\ref{fe}), then
$\widehat{\varepsilon} :F^{^{\mathrm e}}\Rightarrow F$ is the morphism in
$\cgner\hspace{-3pt}^{^\mathrm{u}}\c_n$ given  by the family of 2-cells
$$
\xymatrix@C=9pt@R=9pt{ Fj \ar@/^1pc/[rr]^{\textstyle F^{^{\mathrm e}}_{i,j}}
\ar@/_1pc/[rr]_{\textstyle F_{i,j}} &  \Downarrow \widehat{\varepsilon}_{i,j} & Fi,&0\leq i \leq
j\leq n,}
$$
where $\widehat{\varepsilon}_{i,i}=1_{1_{Fi}}$, $\widehat{\varepsilon}_{i,i+1}=1_{F_{i,i+1}}$, and
then inductively defined by pasting the diagrams
$$
\xymatrix{&&&\ar@{}[dd]|>>>>>>{\textstyle \Downarrow
\widehat{F}_{i,j,j+1}}&\\Fi\ar@{}[rrrru]|-{=}&&Fj\ar[ll]_{\textstyle F^{^{\mathrm
e}}_{i,j}}\ar@/^20pt/[ll]^<<<<<{\textstyle F_{i,j}}&&F(j+1). \ar@/_30pt/[llll]_{\textstyle
F^{^{\mathrm e}}_{i,j+1}}\ar[ll]_{\textstyle F_{j,j+1}}\ar@/^40pt/[llll]^{\textstyle
F_{i,j+1}}\\&\ar@{}[u]|>>>{\textstyle \Downarrow \widehat{\varepsilon}_{i,j}}&&&}
$$

Since $ P_n\widehat{\varepsilon}=1_{P_n}$ and $\widehat{\varepsilon} J_n=1_{J_n}$, it follows that
$P_n:\cgner\hspace{-3pt}^{^\mathrm{u}}\c_n\to \cner\c_n $ is a right adjoint functor to
$J_n:\cner\c_n \to \cgner\hspace{-3pt}^{^\mathrm{u}}\c_n$, with the identity and $\varepsilon$
being the unit and the counit of the adjunction, respectively.

It follows that every $J_n:\cner\c_n \to \cgner\hspace{-3pt}^{^\mathrm{u}}\c_n$, $n\geq 0$, induces
a homotopy equivalence
 on classifying spaces $\class\cner\c_n\to \class \cgner\hspace{-3pt}^{^\mathrm{u}}\c_n$. Therefore, by Theorem
  \ref{p.1.23}, the map induced by the pseudo-simplicial functor
   $\class J:\class\c=\class\cner\c\to \class \cgner\hspace{-3pt}^{^\mathrm{u}}\c$ is a homotopy equivalence as well.

Finally we restrict the above constructions to the simplicial subcategory
$\underline{\mathrm{S}}\c\subseteq \cgner\hspace{-3pt}^{^\mathrm{u}}\c$. Note that, when
 $F:[n]\rightsquigarrow \c$ is a normal homomorphism, that is, an object of $\underline{\mathrm{S}}\c_n$, then
 $\widehat{\varepsilon} :J_nP_n(F)\Rightarrow F$ is an isomorphism. It follows that every
 functor $J_n:\cner\c_n\to \underline{\mathrm{S}}\c_n$ is  an equivalence of categories, with
   $P_n:\underline{\mathrm{S}}\c_n\to \cner\c_n$ being a quasi-inverse for $n\geq 2$ and,
   therefore,
    $J: \cner\c \to \underline{\mathrm{S}}\c$ is actually a pseudo-simplicial equivalence. Hence, the induced map
    $\class J:\class\c\to \class \underline{\mathrm{S}}\c$ is a homotopy equivalence, as claimed.
\end{proof}

Since the triangle  $$\xymatrix{&\cner\c\ar@{~>}[dr]^{\textstyle J}\ar@{~>}[dl]_{\textstyle
J}&\\\cgner\hspace{-3pt}^{^\mathrm{u}}\c&&\underline{\mathrm{S}}\c \ar@{->}[ll]}$$ commutes, we
deduce the following consequence of the above Theorem \ref{psin}, which makes the proof for Theorem
\ref{thit} complete:
\begin{theorem}
For any bicategory $\c$, the simplicial inclusion functor $\underline{\mathrm{S}}\c\hookrightarrow
\cgner\hspace{-3pt}^{^\mathrm{u}}\c$ induces a homotopy equivalence on classifying spaces; thus,
$\class \underline{\mathrm{S}}\c\simeq \class\cgner\hspace{-3pt}^{^\mathrm{u}}\c$.
\end{theorem}

We should stress here that it is not true that the simplicial set of objects of Segal's nerve
 $\underline{\mathrm{S}}\c$ of a bicategory $\c$, that is, the simplicial set
  ${\mathrm S}\c:[p]\mapsto \mathrm{NorHom}([p],\c)$, represents the homotopy type of the bicategory, in contrast
  to what happens with the categorical geometric nerves.

\begin{remark}\label{gp} {\em Let $\Omega^{^{-1}}\hspace{-5pt}\m$ be the delooping bicategory of a monoidal category
$\m=(\m,\otimes,\text{I},\boldsymbol{a},\boldsymbol{l},\boldsymbol{r})$ (see Remark \ref{mono}).
The simplicial space $[n]\mapsto \class (\underline{\mathrm{S}}\Omega^{^{-1}}\hspace{-5pt}\m_n)$
satisfies the hypothesis of Segal's Proposition 1.5 in \cite{segal74}, that is, the space
$\class\underline{\mathrm{S}}\Omega^{^{-1}}\hspace{-5pt}\m_0$ is contractible (since it is a one
point space) and the projection maps
 $p_n=\class P_n:\class (\underline{\mathrm{S}}\Omega^{^{-1}}\hspace{-5pt}\m_n)\to (\class \underline{\mathrm{S}}\Omega^{^{-1}}\hspace{-5pt}\m_1)^n=\class\m^n$
are homotopy equivalences, where $\class\m$ is the classifying space of the underlying category
$\m$ (since every $P_n:\underline{\mathrm{S}}\Omega^{^{-1}}\hspace{-5pt}\m_n\to \m^n$ is an
equivalence of categories, as  observed in the proof of Theorem \ref{psin} above). Therefore, $
\Omega(\class\Omega^{^{-1}}\hspace{-5pt}\m,*)$  is a group completion of $\class \m$.

Note that the same conclusion, with the same discussion as above, can be obtained by using any of
the unitary categorical geometric nerves of the delooping bicategory
$\Omega^{^{-1}}\hspace{-5pt}\m$, that is,
$\cgner\hspace{-3pt}^{^\mathrm{u}}\Omega^{^{-1}}\hspace{-5pt}\m$ or
 $\overline{\nabla}\hspace{-2.5pt}_{^{_\mathrm{u}}}\Omega^{^{-1}}\hspace{-5pt}\m$ instead of the Segal nerve
  for modeling $\class\Omega^{^{-1}}\hspace{-5pt}\m$.
}\qed
\end{remark}

\section{The homotopy colimit theorem for bicategories}
Many properties of the classifying space construction for bicategories, $\c\mapsto \class\c$, are
easily established by using geometric nerves. For example:

\begin{proposition}\label{trans}\begin{enumerate}
\item[(i)] Any lax (resp. oplax) functor between bicategories $F:\b\rightsquigarrow\c$ induces a continuous map $\class F:\class\b\to\class\c$, well defined up to homotopy equivalence.
\item[(ii)] If two lax (resp. oplax) functors between bicategories, $F,F':\b\rightsquigarrow\c$, are related by a lax or an oplax
transformation, $F\Rightarrow F'$, then the induced maps on classifying spaces, $ \class F, \class
F':\class\b\to\class\c$, are homotopic.
\item[(iii)] If a homomorphism of bicategories has a left or right biadjoint, then the induced map on classifying spaces is a homotopy
equivalence. In particular, any biequivalence of bicategories induces a homotopy equivalence on
classifying spaces.
\end{enumerate}
\end{proposition}
\begin{proof} (i) The geometric nerve construction on bicategories $\c\mapsto \gner\c$ (resp. $\c\mapsto \nabla\c$) is functorial on lax (resp. oplax) functors between bicategories. Therefore, Theorem \ref{thit} gives the result.

(ii) Let $F,F':\b\rightsquigarrow\c$ be lax functors and suppose $\alpha:F\Rightarrow F'$ is a lax
transformation. There is a lax functor $H:\b\times [1]\rightsquigarrow \c$ making  the diagram
commutative
\begin{equation}\label{homo}
\xymatrix@C=60pt{\b\times [0]\cong \b \ar[d]_{1\times\delta_0}\ar@{~>}[dr]^{F}& \\ \b\times[1]\ar@{~>}[r]^H&\c,\\
\b\times [0]\cong \b\ar@{~>}[ru]_{F'}\ar[u]^{1\times \delta_1}&}
\end{equation}
that carries a morphism in $\b\times[1]$ of the form $(x\overset{u}\to y,1\to 0):(x,1)\to (y,0)$ to
the composite morphism in $\c$
$$Fx\overset{\textstyle \alpha x}\longrightarrow F'x\overset{\textstyle F'u}\longrightarrow  F'y,$$
and a deformation $(\phi,1_{1\to 0}):(x\overset{u}\to y,1\to 0)\Rightarrow (x\overset{v}\to y,1\to
0)$ to
$$ F'\phi\circ 1_{\alpha x}:
 F'u\circ \alpha x\Rightarrow  F'v\circ\alpha x.$$

For $x\overset{u}\to y\overset{v}\to z$ two composible morphisms in $\b$,  the structure constraint
$$\widehat{H}:H(y\overset{v}\to z,1\to 0)\circ H(x\overset{u}\to y,1\to 1)\Rightarrow H(x\overset{v\circ u}\longrightarrow z,1\to 0)$$
is the deformation obtained by pasting the diagram
$$
\xymatrix@R=50pt{Fx\ar@{}[rd]^<<<<<<<<{\Downarrow\textstyle{
\widehat{\alpha}}}\ar[r]^{\textstyle{Fu}}\ar[d]_{\textstyle{ \alpha x}}&Fy\ar[r]^{\textstyle{
\alpha y}}&F'y\ar[d]^{\textstyle{ F'v}}\\F'x\ar[rru]|{\textstyle {F'u}}\ar[rr]_{\textstyle
F'(v\circ u)}&&F'z,\ar@{}[lu]|>>>>>>>>>>>>>>{\Downarrow\textstyle{ \widehat{F'}}}}
$$
whereas the constraint
$$\widehat{H}:H(y\overset{v}\to z,0\to 0)\circ H(x\overset{u}\to y,1\to 0)\Rightarrow H(x\overset{v\circ u}\longrightarrow  z,1\to 0)$$
is the composite deformation
$$
\xymatrix{ F'v\circ ( F'v \circ \alpha x)\ar@{=>}[r]^{\textstyle\boldsymbol{a}}&(F'v\circ F'u)\circ
\alpha x\ar@{=>}[r]^{\textstyle\widehat{F'}\circ 1}&  F'(v\circ u)\circ \alpha x.}
$$

Applying geometric nerve construction  to diagram (\ref{homo}), we get a diagram of simplicial set
maps
$$
\xymatrix@C=60pt{\gner\b\times \gner[0]\cong \gner\b \ar[d]_{1\times \delta_0}\ar[dr]^{F_*}& \\
\gner\b\times\gner[1]\ar[r]^{H_*}&\gner\c,\\
\gner\b\times \gner[0]\cong \gner\b\ar[ru]_{ F'_*}\ar[u]^{1\times\delta_1}&}
$$
showing that the simplicial maps $F_*,F'_*:\gner\b\to \gner\c$ are made homotopic by $H_*$, whence
the result follows by Theorem \ref{thit}.

The proof is similar for the case in which $\alpha:F\Rightarrow F':\b\rightsquigarrow\c$  is an
oplax transformation, but with a change in the construction of the lax functor
${H:\b\times[1]\rightsquigarrow \c}$ that makes diagram (\ref{homo}) commutative: now define  $H$
such that $H(x\overset{u}\to y, 1\to 0)=\alpha y\circ Fu:Fx\to F'y$.

The statements in (iii) follow from (i) and (ii).\end{proof}

In this section, we mainly want to show how the use of geometric nerves and Theorem \ref{thit} can
be applied to homotopy theory of diagrams of bicategories $\c: I^{^{\mathrm{op}}}\to \bicat $,
through a bicategorical {\em Grothendieck construction}
$$
\xymatrix{\int_I \!\c,}
$$
which assembles all bicategories $\c_i$ as we explain next. Recall that $\bicat$, here, means the
category of bicategories with homomorphisms.

\begin{definition}Let $\c:I^{^{\mathrm{op}}}\to \bicat$ be a functor, $(j\overset{a}\to i)\mapsto   (\c_i\overset{a^*}\to \c_j)$, where $I$ is a small category. Then, the Grothendieck construction on $\c$ is  the bicategory, denoted by $\xymatrix{\int_I\c}$, whose objects are pairs $(x,i)$, where $i$ is an object of $I$ and $x$ is one of the bicategory $\c_i$. A morphism $(u,a):(y,j)\to (x,i)$ in $\int_I\c$ is a pair of morphisms where $a:j\to i$ in $I$ and $u:y\to a^*x$ in $\c_j$. Given two morphisms $(u,a),(u',a'):(y,j)\to(x,i)$, the existence of a  deformation $(u,a)\Rightarrow (u',a')$ requires that $a=a'$,  and then, such a deformation
$$\xymatrix@C=6pt{(y,j)  \ar@/^1.3pc/[rr]^{\textstyle (u,a)} \ar@/_1.3pc/[rr]_{\textstyle (u',a)} & {}_{\textstyle \Downarrow(\phi,a)} &(x,i) }  $$ consists of a deformation  $\xymatrix @C=0.5pc {y  \ar@/^0.8pc/[rr]^{\textstyle u} \ar@/_0.8pc/[rr]_{\textstyle u'} & {}_{\textstyle\Downarrow\phi} &a^*x }$ in $\c_j$. Thus,
$\mbox{Ob}\!\int_I\c=\bigsqcup\limits_{i\in \!\mbox{\scriptsize Ob}I}\mbox{Ob}\c_i$, and the
hom-categories of $\int_I\c$ are
$$\xymatrix{\int_I\c\big((y,j),(x,i)\big)=\bigsqcup\limits_{j\overset{a}\to i}\c_j(y,a^*x),}$$
where the disjoint union is over all arrows $a:j\to i$ in $I$.

For each triplet of objects $(z,k)$, $(y,j)$ and $(x,i)$ of $\int_I\c$, the horizontal composition
functor
$$\bigsqcup\limits_{j\overset{a}\to i}\c_j(y,a^*x) \times \bigsqcup\limits_{k\overset{b}\to j}\c_k(z,b^*y) \overset{\circ}\longrightarrow
\bigsqcup\limits_{k\overset{c}\to i}\c_k(z,c^*i)$$ maps the component at two morphisms $a:j\to i$
and $b:k\to j$ of $I$ into the component at the composite $ab:k\to i$ via the composition
$$\xymatrix@C=16pt{\c_j(y,a^*\!x)\!\times\! \c_k(z,b^*\!y) \ar[r]^-{b^*\times 1} &
\c_k(b^*y,b^*a^*x)\times \c_k(z,b^*y)\ar[r]^-{\circ}&\c_k(z,b^*a^*x), }$$that is,
$$\xymatrix@C=2pt{(z,k)  \ar@/^1.3pc/[rr]^{\textstyle (v,b)} \ar@/_1.3pc/[rr]_{\textstyle (v',b)} & {}_{\textstyle \Downarrow(\psi,b)} &(y,j)  \ar@/^1.3pc/[rr]^{\textstyle (u,a)} \ar@/_1.3pc/[rr]_{\textstyle (u',a)} & {}_{\textstyle \Downarrow(\phi,a)} &(x,i)&\ar@{|->}[rrr]^-{\textstyle \circ} && &&(z,k) \ar@/^1.5pc/[rr]^{\textstyle (b^*u\circ v,ab)} \ar@/_1.5pc/[rr]_{\textstyle (b^*u'\circ v',ab)}&{}_{\textstyle \Downarrow(b^*\phi\circ\psi,ab)}& (x,i).  }  $$

Given any three composible morphisms $(t,\ell)\overset{\textstyle (w,c)}\longrightarrow
(z,k)\overset{\textstyle (v,b)}\longrightarrow  (y,j)\overset{\textstyle (u,a)}\longrightarrow
(x,i)$ in $\int_I\c$, the structure associativity isomorphism $$(u,a)\circ ((v,b)\circ(w,c)) \cong
((u,a)\circ (v,b))\circ (w,c)$$ is provided by pasting, in the bicategory $\c_\ell$, the diagram
$$
\xymatrix{t\ar[rr]^{\textstyle c^*v\circ w}\ar[dd]_{\textstyle w}&\ar@{}[d]|-{\textstyle
\Downarrow\boldsymbol{a}}&c^*b^*y\ar[dd]^{\textstyle c^*b^*u}\\ &\ar@{}[d]|(.65){\textstyle
\Downarrow{\widehat{c}^*}}&\\c^*z\ar@/^2pc/[rr]^{\textstyle c^*b^*u\circ c^*v}\ar[rr]_{\textstyle
c^*(b^*u\circ v)}&&c^*b^*a^*x,}
$$
where $\boldsymbol{a}$ is the associativity isomorphism in $\c_\ell$ and $\widehat{c}^*$ is the
structure isomorphism of the homomorphism $c^*:\c_k\rightsquigarrow\c_\ell$.

The identity morphism for each object $(x,i)$ in $\int_I\c$ is $$(1_x,1_i):(x,i)\to(x,i),$$ and,
for each morphism $(u,a):(y,j)\to (x,i)$, the left unit constraint
$$(1_x,1_i)\circ(u,a)=(a^*1_x\circ u,a)\cong (u,a)$$ is provided by the left unit  constraint of $\c_j$ and the unit structure constraint of the homomorphism $a^*:\c_i\rightsquigarrow\c_j$ by pasting
$$
\xymatrix@R=18pt{ &&a^*x\ar[rrdd]^{\textstyle a^*1_x}\ar@/_1.7pc/[rrdd]|(0.3){\textstyle
1_{a^*x}}&&\\&&&&\\y\ar@{}[rrrru]|(0.7){\textstyle
\Downarrow\widehat{a}^*}\ar@{}[rru]_(0.8){\textstyle \boldsymbol{l}\Downarrow}\ar[uurr]^{\textstyle
u}\ar[rrrr]_{\textstyle u}&&&&a^*x,}
$$
whereas the right unit constraint
$$(u,a)\circ(1_y,1_j)=(u\circ 1_y,a)\cong (u,a)$$
is directly given by the right unit constraint $\boldsymbol{r}:u\circ 1_y\Rightarrow u$ of $\c_j$
at $u$.
\end{definition}

The well-known Homotopy Colimit Theorem by Thomason, that is, Theorem \ref{thoma}, admits the
following generalization  to diagrams of bicategories $\c:I^{^{\mathrm{op}}}\to \bicat$ in Theorem
\ref{hct}, where it is established how the classifying space of the bicategory Grothendieck
construction  $\class\int_I\c$ can be thought of as a homotopy colimit of the spaces $\class\c_i$
that arise from the initial input data $i\mapsto \c_i$ given by the functor $\c$.

\begin{theorem}\label{hct}
Suppose a category $I$ is given. To every functor $\c:I^{^{\mathrm{op}}}\to\bicat$ there exists a
natural weak homotopy equivalence of simplicial sets
\begin{equation}\label{eta} \xymatrix{\eta: \hoco_I\gner\c\longrightarrow
\gner\int_I\!\c,}\end{equation}  where $\gner\c:I^{^{\mathrm{op}}}\to\sset$ is the diagram of
simplicial sets obtained by composing $\c$ with the geometric nerve functor $\gner:\bicat\to
\sset$.
\end{theorem}
\begin{proof} To show an explicit description of the weak equivalence $\eta$ in the theorem, we shall first  explicitly  describe the simplicial sets $\hoco_I\gner\c$ and $\gner\int_I\c$.

On one hand, the simplicial set $\hoco_I\gner\c$ is the diagonal of the bisimplicial set
\begin{equation}\label{s}
S=\bigsqcup_{G\in\gner I}\gner\c_{G 0}=\bigsqcup_{G:[q]\to I}\lfunc([p],\c_{G 0}),
\end{equation}
whose $(p,q)$-simplices are pairs
$$
(F,G)
$$
consisting of a functor $G:[q]\to I$ and a lax functor $F:[p]\rightsquigarrow \c_{G0}$. If
$\alpha:[p']\to [p]$ and $\beta:[q']\to [q]$ are maps in the simplicial category, then the
respective horizontal and vertical induced maps $\alpha^{*_h}: S_{p,q}\to S_{p',q}$ and
$\beta^{*_v}: S_{p,q}\to
S_{p,q'}$ are defined by $$\left\{\begin{array}{l}\alpha^{*_h}(F,G)=(F\alpha,G),\\[6pt] \beta^{*_v}(F,\sigma)=
(G_{0,\beta 0}^*F,G\beta),\end{array}\right.$$ where $G_{0,\beta 0}^*F:[p]\rightsquigarrow
\c_{G\beta 0}$ is the lax functor obtained by the composition of $F$ with the homomorphism of
bicategories $G_{0,\beta 0}^*:\c_{G0}\rightsquigarrow \c_{G\beta0}$ attached in diagram
$\c:I^{{^\mathrm{op}}}\to\bicat $ at the morphism $G_{0,\beta0}:G\beta0\to G0$ of $I$. In
particular, the horizontal face maps are given by $$d_i^h(F,G)=(Fd^i, G), \hspace{0.2cm} \text{ for
} 0\leq i\leq p,$$ and the vertical ones by
$$d_j^v(F,G)=(F, Gd^j), \hspace{0.2cm} \text{ for } 1\leq j\leq q,$$
while
$$d_0^v(F,G)=(G_{0,1}^*F, Gd^0).$$

On the other hand, a $p$-simplex of $\gner\int_I\c$ is a  lax functor $[p]\rightsquigarrow
\int_I\c$, which can be described as a pair
$$(F',G),$$ where $G:[p]\to I$ is a functor,
that is, a $p$-simplex of $\gner I$, and $F':[p]\rightsquigarrow \c$ is a {\em $G$-crossed lax
functor} (cf. \cite[\S 4.1]{cegarra3}), that is, a family
\begin{equation}\label{xlf}
F'=\{F'\!i,F'_{i,j},\widehat{F}'_{i,j,k},\widehat{F}'_{i}\}_{0\leq i\leq j\leq k\leq p}
\end{equation}
in which each $F'\!i$ is an object of the bicategory $\c_{Gi}$, each $F'_{i,j}:F'\!j\rightarrow
G_{i,j}^*F'\!i$ is a morphism in $\c_{Gj}$, each
$\widehat{F}'_{i,j,k}:G_{\hspace{-3pt}j,k}^*F'_{i,j}\circ F'_{j,k}\Rightarrow F'_{i,k}$ is a
deformation in $\c_{Gk}$
$$
\xymatrix{\ar@{}[drr]|(.6)*+{\textstyle \Downarrow \widehat{F}'_{i,j,k}}
                & G_{\hspace{-3pt}j,k}^*F'\!j \ar[dl]_{\textstyle G_{\hspace{-3pt}j,k}^*F'_{i,j}}             \\
 G_{\hspace{-3pt}j,k}^*G_{\hspace{-2pt}i,j}^*F'\!i=G_{\hspace{-2pt}i,k}^*F'\!i  & &     F'\!k\,,   \ar[ul]_{\textstyle F'_{j,k}} \ar[ll]^{\textstyle F'_{i,k}}}
$$
and each $\widehat{F}'_{i}:1_{F'\!i}\Rightarrow F'_{i,i}$ is a deformation in $\c_{Gi}$. These data
are required to satisfy that, for $0\leq i\leq j\leq k\leq \ell \leq p$, the diagram of
deformations

\begin{equation}\label{ud}
\xymatrix@C=37pt{ (G_{\hspace{-3pt}j,\ell}^*F'_{i,j}\!\circ\! G_{\hspace{-3pt}k,\ell}^*F'_{j,k})\!\circ\! F'_{k,\ell}\ar@{=>}[r]^-{\textstyle \widehat{G}_{\hspace{-3pt}k,\ell}^*\!\circ\! 1}&G^*_{\hspace{-3pt}k,\ell}(G^*_{\hspace{-3pt}j,k}F'_{i,j}\!\circ\! F'_{j,k})\!\circ\! F'_{k,\ell}\ar@{=>}[r]^-{\textstyle G^*_{\hspace{-3pt}k,\ell}\widehat{F}'_{i,j,k}\!\circ\! 1}&G^*_{k,\ell}F'_{i,k}\!\circ\! F'_{k,\ell}\ar@{=>}[d]^{\textstyle \widehat{F}'_{i,k,\ell}}\\
G_{\hspace{-3pt}j,\ell}^*F'_{i,j}\!\circ\!(G_{\hspace{-3pt}k,\ell}^*F'_{j,k}\!\circ\!
F'_{k,\ell})\ar@{=>}[u]^-{\textstyle\boldsymbol{a}} \ar@{=>}[r]^-{\textstyle
1\!\circ\!\widehat{F}'_{j,k,\ell}}&G^*_{\hspace{-3pt}j,\ell}F'_{i,j}\circ \!F'_{j,\ell}
\ar@{=>}[r]^-{\textstyle \widehat{F}'_{i,j,\ell}}&F'_{i,\ell} }
\end{equation}
commutes in the bicategory $\c_{G\ell}$, and, for any $0\leq i\leq j\leq p$, both diagrams below
commute in $\c_{Gj}$.
\begin{equation}\label{ud2}
\xymatrix@C=40pt{ G^*_{\hspace{-2pt}i,j}F'_{i,i}\!\circ\!F'_{i,j}\ar@{=>}[d]_-{\textstyle
\widehat{F}'_{i,i,j}}&G^*_{\hspace{-2pt}i,j}1_{F'\hspace{-1pt}i}\!\circ\!F'_{i,j}
\ar@{=>}[l]_-{\textstyle G^*_{\hspace{-1pt}i,j}\widehat{F}'_{i}\!\circ\! 1}\\
F'_{i,j}&1_{G^*_{\hspace{-1pt}i,j}F'\!i}\!\circ\!F'_{i,j}\ar@{=>}[u]_{\textstyle
\widehat{G}^*_{\hspace{-1pt}i,j}\!\circ\!1}\ar@{=>}[l]_{\textstyle \boldsymbol{r}} }\hspace{0.6cm}
\xymatrix@C=15pt{F'_{i,j}\!\circ\!F'_{j,j}\ar@{=>}[dr]_{\textstyle
\widehat{F}'_{i,j,j}}&&F'_{i,j}\!\circ\!1_{F'\!j}\ar@{=>}[ll]_-{\textstyle
1\!\circ\!\widehat{F}'_j} \ar@{=>}[dl]^{\textstyle \boldsymbol{l}}\\&F'_{i,j}&}
\end{equation}

If $\alpha:[p']\to [p]$ is a map in the simplicial category, then the induced map
$\alpha^*:(\gner\int_I\!\c)_p\to (\gner\int_I\!\c)_{p'}$ associates to a $p$-simplex $(F',G)$, as
above, the $p'$-simplex
$$\alpha^*(F',G)=(F'\alpha,G\alpha),$$
in which $F'\alpha=\{F'\!{\alpha i},F'_{\alpha i,\alpha
j},\widehat{F}'_{\alpha_i,\alpha_j,\alpha_k}, \widehat{F}'_{\alpha i}\}_{0\leq i\leq j\leq k\leq
p'}$.

The simplicial map (\ref{eta})  $$\xymatrix{\eta: \hoco_I\gner\c\longrightarrow \gner\int_I\!\c,}$$
is then given on a $p$-simplex $(F,G)$ of $\hoco_I\gner\c$ by $$\eta(F,G)=(F',G),$$ where
$F':[p]\rightsquigarrow \c$ is the $G$-crossed lax functor, as in (\ref{xlf}), defined by the
objects $F'\!i=G^*_{0,i}Fi$, the morphisms $F'_{i,j}=G^*_{0,j}F_{i,j}:F'\!j\rightarrow
G_{i,j}^*F'\!i$, and the deformations  $\widehat{F}'_{i,j,k}:G_{\hspace{-3pt}j,k}^*F'_{i,j}\circ
F'_{j,k}\Rightarrow F'_{i,k}$ and $\widehat{F}'_{i}:1_{F'\!i}\Rightarrow F'_{i,i}$, which are,
respectively, the composites
$$\begin{array}{l}
\xymatrix@C=42pt{\widehat{F}'_{i,j,k}:\ G^*_{0,k}F_{i,j}\circ G^*_{0,k}F_{j,k}\ar@{=>}[r]^-{\textstyle \widehat{G}_{0,k}^*}&G^*_{0,k}(F_{i,j}\circ F_{j,k})\ar@{=>}[r]^-{\textstyle G_{0,k}^*\widehat{F}_{i,j,k}}&G^*_{0,k}F_{i,k},}\\[4pt]
\xymatrix@C=35pt{\widehat{F}'_{i}: \ 1_{G^*_{0,i}Fi}\ar@{=>}[r]^-{\textstyle
\widehat{G}^*_{0,i}}&G^*_{0,i}1_{Fi}\ar@{=>}[r]^-{\textstyle
G^*_{0,i}\widehat{F}_i}&G_{0,i}^*F_{i,i}.}
\end{array}
$$

To prove that this map $\eta$ is a weak equivalence, our strategy now is to apply the weak homotopy
equivalences (\ref{Phi}) on the bisimplicial set  $S$, defined in $(\ref{s})$. Since $\diag S=
\hoco_I\gner\c$, we have a weak homotopy equivalence ${\Phi: \hoco_I\gner\c \to \overline{W}S}$,
and the proof will be complete once we show a simplicial isomorphism
$$\xymatrix{\Psi: \overline{W}S\ \cong \ \gner\int_I\!\c}$$ making  the diagram of simplicial sets commutative:
\begin{equation}\label{tf}\xymatrix{\hoco_I\gner\c\ar[rr]^{\textstyle \eta}\ar[rd]_{\textstyle \Phi}^\simeq&&\gner\int_I\!\c\\&\overline{W}S.\ar[ru]_{\textstyle \Psi}^\cong&}
\end{equation}

For, note that a $p$-simplex of $\overline{W} S$, say $\chi$, can be described as a list of pairs
$$
\chi=\big(\big(F^{^{(0}}\hspace{-3pt},G^{^{(p}}\big), \dots,
\big(F^{^{(m}}\hspace{-3pt},G^{^{(p-m}}\big),\dots,\big(F^{^{(p}}\hspace{-3pt},G^{^{(0}}\big)\big),
$$
in which each $G^{^{(p\text{-}m}}:[p-m]\to I$ is a functor and each $F^{^{(m}}:[m]\rightsquigarrow
\c_{G^{^{(p\text{-}m}}\!0}$ is a lax functor, such that the equalities
$$\begin{array}{l}G^{^{(p\text{-}m}}d^0=G^{^{(p\text{-}m\text{-}1}},\hspace{0.6cm}
F^{^{(m\text{+}1}} d^{m+1}={G^{^{(p\text{-}m}}_{0,1}}^*F^{^{(m}},\end{array}$$ hold for all $0\leq
m<p$. Writing $G^{^{(p}}:[p]\to I$ simply as $G:[p]\to I$, an iterated use of the above equalities
proves that
$$
G^{^{(p\text{-}m}}=G d^{0}\overset{(m}\cdots d^0 :[p-m]\to I,
$$
for $0\leq m\leq p$, and
$$
F^{^{(m\text{+}1}} d^{m+1}\cdots d^{k+1}=G_{\hspace{-2pt}k,m+1}^*F^{^{(k}}:[k]\rightsquigarrow \c_{G(m+1)},$$ for $0\leq k\leq m<p$. These latter equations mean that  $$ \left\{\begin{array}{ll}F^{^{(j}}\!i=G_{\hspace{-2pt}i,j}^*\,F^{^{(i}}\!i& \text{for } \ i\leq j,\\[4pt]
F^{^{(k}}_{i,j}=G_{\hspace{-2pt}j,k}^*\,F^{^{(j}}_{i,j}& \text{for } \ i\leq j\leq k,\\[4pt]
\widehat{F}^{^{(\ell}}_{i,j,k}=G_{\hspace{-2pt}k,\ell}^*\ \widehat{F}^{^{(k}}_{i,j,k}& \text{for }\  i\leq j\leq k\leq \ell,\\[4pt]
\widehat{F}^{^{(j}}_i=G_{\hspace{-3pt}i,j}^* \widehat{F}^{^{(i}}_i& \text{for } i\leq j,
\end{array}\right.
$$ whence we see how the $p$-simplex $\chi$ of $\overline{W} S$ is uniquely determined by $G:[p]\to I$, the objects
$F^{^{(i}}\!i$ of $\c_{Gi}$,  the morphisms $F^{^{(j}}_{i,j}:F^{^{(j}}\!j\to
F^{^{(j}}\!i=G_{\hspace{-2pt}i,j}^*F^{^{(i}}\!i$ of $\c_{Gj}$, and the deformations
$\widehat{F}^{^{(k}}_{i,j,k}:F^{^{(k}}_{i,j}\circ
F^{^{(k}}_{j,k}=G_{\hspace{-2pt}j,k}^*F^{^{(j}}_{i,j}\circ F^{^{(k}}_{j,k}\Rightarrow
F^{^{(k}}_{i,k}$ in $\c_{Gk}$, and $\widehat{F}^{^{(i}}_i:1_{F^{^{(i}}\!i}\Rightarrow
F^{^{(i}}_{i,i}$ in $\c_{Gi}$,  all for $0\leq i\leq j\leq k\leq p$. At this point, we observe that
there is a normal $G$-crossed lax functor
$F'=\{F'\!i,F'_{i,j},\widehat{F}'_{i,j,k},\widehat{F}_i'\}:[p]\rightsquigarrow \c$, as in
(\ref{xlf}), defined just by putting $F'\!i=F^{^{(i}}\!i$, $F'_{i,j}=F^{^{(j}}_{i,j}$,
$\widehat{F}'_{i,j,k}=F^{^{(k}}_{i,j,k}$ and $\widehat{F}'_i=F^{^{i}}_i$ (the commutativity of
diagrams (\ref{ud}) and (\ref{ud2}) follows from $F^{^{(\ell}}$ and $F^{^{(j}}$ being  lax
functors). Thus, the $p$-simplex $\chi\in \overline{W}S$ defines the $p$-simplex  $(F',G)$ of
$\gner\int_I\!\c$, which itself uniquely determines  $\chi$.

In this way, we get an injective simplicial map $$\xymatrix{\Psi:\overline{W}S\to
\gner\int_I\!\c,}$$
$$\big(\big(F^{^{(0}}\hspace{-3pt},G^{^{(p}}\big), \dots,\big(F^{^{(p}}\hspace{-3pt},G^{^{(0}}\big)\big)\mapsto (F',G)=
\big(\{F^{^{(i}}\!i,F^{^{(j}}_{i,j}, \widehat{F}^{^{(k}}_{i,j,k}, \widehat{F}^{^{(i}}_i\},
G^{^{(p}}\big)$$ which is also surjective, that is, actually an isomorphism, as we can see by
retracing our steps:

 To any pair $(F',G)$ describing a $p$-simplex of $\gner\int_I\!\c$, that is, with $G:[p]\to I$ a functor and $F'=\{F'\!i,F'_{i,j},\widehat{F}'_{i,j,k},\widehat{F}_i'\}
 :[p]\rightsquigarrow\c$ a $G$-crossed lax
 functor, we associate the $p$-simplex
 $\chi=\big(\big(F^{^{(m}}\hspace{-2pt},G^{^{(p\text{-}m)}}\big)\big)$ of $\overline{W}S$, where, for each $0\leq m\leq p$,
 $G^{^{(p\text{-}m}}:[p-m]\to I$ is the
 composite $[p-m]\overset{(d^0)^{m}}\to [p] \overset{ G}\to I$,
and the lax functor $F^{^{(m}}:[m]\rightsquigarrow \c_{G^{^{(p\text{-}m}}\!0}$ is defined by the
objects $F^{^{(m}}\!i=G_{\hspace{-2pt}i,m}^*F'\!i$, the morphisms
$F^{^{(m}}_{i,j}=G_{\hspace{-2pt}j,m}^*F'_{i,j}:F^{^{(m}}\!{j}\to F^{^{(m}}\!{i}$, and the
deformations $\widehat{F}^{^{(m}}_{i,j,k}:F^{^{(m}}_{i,j}\circ F^{^{(m}}_{j,k}\Rightarrow
F^{^{(m}}_{i,k}$ and $F^{^{(m}}_{i}:1_{F^{^{(m}}\!i}\Rightarrow F^{^{(m}}_{i,i}$ which are
respectively given as the compositions
$$\begin{array}{l}
\xymatrix@C=42pt{\widehat{F}^{^{(m}}_{i,j,k}:\ G^*_{\hspace{-2pt}j,m}F'_{i,j}\circ G^*_{\hspace{-2pt}k,m}F'_{j,k}\ar@{=>}[r]^-{\textstyle \widehat{G}_{\hspace{-2pt}k,m}^*}&G^*_{\hspace{-2pt}k,m}(G^*_{\hspace{-2pt}j,k}F'_{i,j}\circ F'_{j,k})\ar@{=>}[r]^-{\textstyle G_{k,m}^*\widehat{F}'_{i,j,k}}&G^*_{\hspace{-2pt}k,m}F'_{i,k},}\\[4pt]
\xymatrix@C=35pt{\widehat{F}^{^{(m}}_{i}: \ 1_{G^*_{i,m}F'\!i}\ar@{=>}[r]^-{\textstyle
\widehat{G}^*_{\hspace{-2pt}i,m}}&G^*_{\hspace{-2pt}i,m}1_{F'\!i}\ar@{=>}[r]^-{\textstyle
G^*_{\hspace{-2pt}i,m}\widehat{F}'_i}&G_{\hspace{-2pt}i,m}^*F'_{i,i}.}
\end{array}
$$
One easily checks that $\Psi(\chi)=(F',G)$, whence we conclude that the simplicial map $\Psi$ is an
isomorphism as claimed.

Finally, since the map $\eta$ in the theorem occurs in the commutative diagram (\ref{tf}), that is,
$\eta=\Psi\,\Phi$, where $\Psi$ is an isomorphism and $\Phi$ a weak homotopy equivalence, the proof
is complete.
\end{proof}

For any functor $\c: I^{^\mathrm {op}} \to \bicat $, the bicategory $\int_I\!\c$ assembles all
bicategories $\c_i$ in the following precise sense: There is a projection 2-functor
$$
\xymatrix{\pi :\int_I\!\c \to I,}
$$
given by
$$\xymatrix@C=6pt{(y,j)  \ar@/^1.3pc/[rr]^{\textstyle (u,a)} \ar@/_1.3pc/[rr]_{\textstyle (v,a)} & {}_{\textstyle \Downarrow(\phi,a)} &(x,i)&\overset{\textstyle \pi}\mapsto\underset{~}~& j\ar[rr]^{a}&&i \,,}  $$
and, for each object $j$ of $I$, there is a pullback square of bicategories
\begin{equation}\label{square}
\xymatrix{\c_j\ar[r]^{\textstyle \iota}\ar[d]&\int_I\!\c\ar[d]^{\textstyle \pi}\\ [0]
\ar[r]^{\textstyle j}&I}
\end{equation}
where  $\iota:\c_j\to \int_I\!\c$ is the embedding homomorphism defined by
$$\xymatrix@C=3pt{y \ar@/^0.7pc/[rr]^{\textstyle u} \ar@/_0.7pc/[rr]_{\textstyle v} & {}_{\textstyle \Downarrow\phi} &x}\ \overset{\textstyle \iota}\mapsto  \
\xymatrix@C=4pt{(y,j)  \ar@/^1.3pc/[rr]^{\textstyle (u,1_j)} \ar@/_1.3pc/[rr]_{\textstyle (v,1_j)}
& {}_{\textstyle \Downarrow(\phi,1_j)} &(x,j)} .$$ Thus, $\c_j\cong \pi^{-1}(j)$, the fibre
bicategory of $\pi$ at $j$.

After Quillen's Lemma  \cite[pag. 90]{quillen}, \cite[\S IV, Lemma 5.7]{g-j}, the following result
is a consequence of our Theorem \ref{hct}:

\begin{theorem}\label{hc2} Suppose that $\c:I^{^{\mathrm{op}}}\to \bicat$ is a diagram of bicategories such that the induced map
$\class a^*:\class \c_{i}\to \class \c_{j}$, for each morphism $a:j\to i$ in  $I$,  is a homotopy
equivalence. Then, for every object $j$ of $I$, the square induced by $(\ref{square})$ on
classifying spaces
$$
\xymatrix{\class\c_j\ar[r]\ar[d]&\class \int_I\!\c\ar[d]\\ \star \ar[r]^j&\class I}
$$
is homotopy cartesian. Therefore, for each object $y\in \c_j$  there is an induced long exact
sequence on homotopy groups, relative to the base points $y$ of $\class\c_j$, $(y,j)$ of $\class
\!\int_I\!\c$, and $j$ of $\class I$,
$$\xymatrix{\cdots \to \pi_{n+1}\class I\to\pi_n\class\c_j\to\pi_n\class \!\int_I\!\c\to\pi_n\class I\to\cdots.}$$
\end{theorem}

\begin{remark}{\em  Let $\m: I^{^\mathrm {op}} \to \mathbf{MonCat}$, $(j\overset{a}\to i)\mapsto   (\m_i\overset{a^*}\to \m_j)$, be a diagram of monoidal categories. It follows from Theorems \ref{thit} and \ref{hct} that the  homotopy type of $\m$ is modeled by the bicategory $$\xymatrix{\int_I\!\Omega^{^{-1}}\hspace{-5pt}\m,}$$
where $\Omega^{^{-1}}\hspace{-5pt}:\mathbf{MonCat}\to \bicat$ is the delooping embedding (see
Remarks \ref{mono} and \ref{monf}).

Notice that $\xymatrix{\int_I\!\Omega^{^{-1}}\hspace{-5pt}\m}$ is a genuine bicategory: It has the
same objects $i$ as the category $I$,  its hom-categories are
$$\xymatrix{\int_I\!\Omega^{^{-1}}\hspace{-5pt}\m(j,i)=\bigsqcup\limits_{j\overset{a}\to i}\m_j=\m_j\times I(j,i),}$$
where $\m_j$ is the underlying category of the monoidal category equally denoted, and its
horizontal compositions are given by
$$\xymatrix@C=2pt{k  \ar@/^1.3pc/[rr]^{\textstyle (v,b)} \ar@/_1.3pc/[rr]_{\textstyle (v',b)} & {}_{\textstyle \Downarrow(\psi,b)} &j  \ar@/^1.3pc/[rr]^{\textstyle (u,a)} \ar@/_1.3pc/[rr]_{\textstyle (u',a)} & {}_{\textstyle \Downarrow(\phi,a)} &i&\ar@{|->}[rrr]^-{\textstyle \circ} && &&k \ar@/^1.5pc/[rr]^{\textstyle (b^*u\otimes v,ab)} \ar@/_1.5pc/[rr]_{\textstyle (b^*u'\otimes v',ab)}&{}_{\textstyle \Downarrow(b^*\varphi\otimes\psi,ab)}& i.  }  $$
Hence, the reader interested in the study of classifying spaces of monoidal categories can find in
the above fact a good reason to also be interested in the study of classifying spaces of
bicategories.

\vspace{0.2cm} Recall that a notion of nerve $ \mbox{Ner}_I\m$, for a diagram of monoidal
categories   $\m: I^{^\mathrm {op}} \to \mathbf{MonCat}$, was defined in \cite[(66)]{cegarra3}  as
follows:  A `2-cocycle' of $I$ with coefficients in $\m$ is  a system of data $(Y,f)$ consisting
of:

- For each arrow $j\overset{a}\to i$ in $I$, an object $Y_a\in \m_j$.

- For each pair of composible arrows in  $I$, $k\overset{b}\to j\overset{a}\rightarrow i$, a
morphism in $\m_k$
$$\xymatrix{b^*Y_a\otimes Y_b\ar[r]^-{\textstyle f_{a,b}}&Y_{ab},}$$
such that  $Y_{1_j}=\mathrm{I}$ (the unit object of $\m_j$), the morphisms
$f_{1,a}:a^*\mathrm{I}\otimes Y_a\to Y_a$ and $f_{a,1}:Y_a\otimes \mathrm{I} \to Y_a$ are the
canonical isomorphisms given by the unit constrains of the monoidal category $\m_j$ and the
monoidal functor $a^*$, and for any three composable triplet, $\ell\overset{c}\to k \overset{b}\to
j\overset{a}\to i$, of morphisms in $I$, the diagram in $\m_\ell$
$$\xymatrix@C=50pt@R=25pt{(c^*b^*Y_a\otimes c^*Y_b)\otimes Y_c\ar[r]^{\cong}&c^*(b^*Y_a\otimes Y_b)\otimes Y_c\ar[r]^-{\textstyle c^*\!f_{a,b}\otimes 1}&c^*Y_{ab}\otimes Y_c \ar[d]^{\textstyle f_{ab,c}}\\c^*b^*Y_a\otimes(c^*Y_b\otimes Y_c)\ar[u]^{\cong}\ar[r]^-{\textstyle 1\otimes f_{b,c}}&c^*b^*Y_a\otimes Y_{bc}\ar[r]^-{\textstyle f_{a,bc}}& Y_{abc}}
$$
(where the unnamed isomorphisms are canonical) is commutative. Then,  $\mbox{Ner}_I\m$, the nerve
of the diagram,  is defined as the simplicial set
\begin{equation}\label{4.6} \mbox{Ner}_I\m:\ [n]\mapsto \bigsqcup_{G:[n]\to I} Z^2\big([n],\m\,G\big)\,,
\end{equation}
where $G:[n]\to I$ is any functor and $Z^2\big([n], \m\,G\big)$ is the set of $2$-cocycles of $[n]$
in the composite functor $[n]\overset{G}\to I\overset{\m}\to \mathbf{MonCat}$.

It is a fact that this simplicial set (\ref{4.6}) actually represents the homotopy type of the
diagram of monoidal categories $\m$ since a straightforward comparison shows the existence of a
natural isomorphism of simplicial sets $$\xymatrix{\mbox{Ner}_I\m\cong
\Delta^{\hspace{-2pt}^\mathrm{u}}\!\int_I\!\Omega^{^{-1}}\hspace{-5pt}\m,}$$ between
$\mbox{Ner}_I\m$ and the normal geometric nerve (\ref{ngn}) of the bicategory
$\int_I\!\Omega^{^{-1}}\hspace{-5pt}\m$. }\qed
\end{remark}


\begin{thebibliography}{99}
\bibitem{artin-mazur} {\bf M Artin, B Mazur,} \emph{ On the Van Kampen Theorem,} Topology 5 (1966)
179-189\,\,\,\, MR0192495

\bibitem{benabou} {\bf J B\'{e}nabou,} \emph{ Introduction to bicategories} in ``Reports of the
Midwest Category Seminar" Lecture Notes in Math. 47, Springer-Verlag, Berlin-New York (1967) 1-77
\,\,\,\, MR0220789

\bibitem{bousfield-kan} {\bf A K Bousfield, D M   Kan,} \emph{ Homotopy limits, completions and localizations,}
 Lecture Notes in Math. 304, Springer-Verlag, Berlin-New York  (1972) v+348 pp. \,\,\,\,MR0365573

\bibitem{bre92} {\bf L Breen,} \emph{ Th\'{e}orie de Schreier sup\'{e}rieure,}  Ann. Scient. Ec. Norm. Sup.
 4$^e$ s\'{e}rie 25 (1992) 465-514\,\,\,\, MR1191733

\bibitem{b-c}{\bf M  Bullejos, A M  Cegarra,}  \emph{ On the geometry of 2-categories and their classifying spaces,}
 $K$-Theory (3) 29 (2003) 211-229 \,\,\,\, MR2028502

\bibitem{b-c2} {\bf M Bullejos, A M  Cegarra,}  \emph{ Classifyng spaces for monoidal categories through Geometric
 Nerves,} Canad. Math. Bull. (3) 47 (2004) 321-331 \,\,\,\, MR2072592

\bibitem{Ca-Ceg} {\bf P Carrasco, A M  Cegarra,}  \emph{ (Braided) Tensor
structures on homotopy groupoids and nerves of (braided) categorical groups,} Comm. in Algebra
 24 (1996) 3995-4058 \,\,\,\, MR1414569

\bibitem{ceg-gar} {\bf A M   Cegarra, A R  Garz\'{o}n,}  \emph{ Homotopy classification of categorical torsors,}
 Appl. Cat. Structures 9 (2001) 465-496  \,\,\,\,MR1865612

\bibitem{cegarra1} {\bf A M Cegarra, J  Remedios,}  \emph{ The relationship between the diagonal
 and the bar constructions on a bisimplicial set,} Topology Appl. (1) 153 (2005) 21-51 \,\,\,\, MR2172033



\bibitem{cegarra2} {\bf A M Cegarra, J   Remedios,}  \emph{The behaviour of the $\overline W$-construction
on the homotopy theory of bisimplicial sets,} Manuscripta Math. (4) 124 (2007) 427-457 \,\,\,\,
MR2357792

\bibitem{cegarra3} {\bf A M Cegarra, E  Khadmaladze,}  \emph{ Homotopy classification of graded Picard
 categories,} Adv. Math. (2) 213 (2007) 644-686 \,\,\,\, MR2332605

\bibitem{duskin} {\bf J Duskin,}  \emph{Simplicial matrices and the nerves of weak $n$-categories. I.
 Nerves of bicategories,}  Theory Appl. Cat. (10) 9 (2002) 198-308 \,\,\,\, MR1897816


\bibitem{fiedorowicz} {\bf Z Fiedorowicz,}  \emph{Classifying spaces of topological monoids and categories,}
 Amer. J. Math. (2) 106 (1984) 301-350 \,\,\,\, MR0737777

\bibitem{g-j} {\bf P G Goerss, J F Jardine,} \emph{Simplicial homotopy theory,} Progress in Mathematics
 174 Birkh\"{a}user Verlag, Basel (1999) xvi+510 pp. \,\,\,\,MR1711612


\bibitem{g-g} {\bf R Garner, N  Gurski,}  \emph{ The low-dimensional structures that tricategories form,}
To appear in Mathematical Proc. Camb. Phil. Soc.  ( arXiv:0711.1761)


\bibitem{grothendieck} {\bf A Grothendieck,}  \emph{ Cat\'{e}gories fibr\'{e}es et d\'{e}scente,} SGA I
expos\'{e} VI, Lecture Notes in Math. 224 Springer, Berlin (1971) 145-194

\bibitem{gurski} {\bf N Gurski,}  \emph{ An algebraic theory of tricategory,} Ph.D. thesis, March 9, 2007


\bibitem{gurski2} {\bf N Gurski}  \emph{ Nerves of bicategories as estratified simplicial sets,}
 J. Pure Appl. Algebra (2008), In Press, Corrected Proof available online  doi:10.1016/j.jpaa.2008.10.011.



\bibitem{g-p-s} {\bf R Gordon, A J   Power, R Street,}  \emph{ Coherence for tricategories,} Mem. Amer.
Math. Soc. (558) 117 (1995) \,\,\,\, MR1261589
\bibitem{hinich} {\bf V A Hinich, V V  Schechtman,}  \emph{ Geometry of a category of complexes and algebraic
K-theory,} Duke Math. Journal 52 (1985) 339-430 \,\,\,\, MR0792180

\bibitem{jardine} {\bf J F Jardine,}  \emph{ Supercoherence,} J. Pure
Appl. Algebra 75  (1991) 103-194 \,\,\,\, MR1138365

\bibitem{joyal}  {\bf  A Joyal, R Street}  \emph{ Braided Tensor Categories,} Adv. Math. 102 (1993)
20-78 \,\,\,\, MR1250465

\bibitem{k-v} {\bf M Kapranov, V  Voevodsky,} \emph{ 2-Categories and Zamolodchikov tetrahedra equations,} in
Proc. Symp. Pure Math.  56 p.2  AMS, Providence, (1994) 177-260 \,\,\,\, MR1278735

\bibitem{lack} {\bf S Lack,} ICONS (arXiv: 0711.4657).

\bibitem{lack-paoli} {\bf S Lack, S.  Paoli,}  \emph{ 2-nerves for bicategories,} K-Theory 38 (2008)
153-175 \,\,\,\, MR2366560

\bibitem{quillen} {\bf D Quillen,}  \emph{ Higher algebraic K-theory:I,} in Algebraic K-theory I, Lecture
Notes in Math. 341 Springer-Verlag  (1973) 77-139 \,\,\,\, MR0338129

\bibitem{maclane} {\bf S Mac Lane}   \emph{Categories for the working mathematician,} GTM  5 2nd Edition,
Springer (1998) \,\,\,\, MR1712872


\bibitem{may67} {\bf J P  May,}  \emph{ Simplicial Objects in Algebraic Topology,} Van Nostrand, Princeton, New York,
1967 \,\,\,\, MR0222892

\bibitem{may72} {\bf J P May,}  \emph{ The Geometry of Iterated Loop Spaces,}  Lecture Notes in Math.  271
(1972) \,\,\,\, MR0420610

\bibitem{may} {\bf J P May,}  \emph{ Pairing of categories and spectra,}  J.
 Pure Appl. Algebra  19 (1980) 299-346  \,\,\,\, MR0593258

\bibitem{moerdijk-svensson} {\bf I Moerdijk, J A  Svensson,}  \emph{Algebraic classification of equivariant
homotopy 2-types,I,} J. Pure Appl. Algebra  89 (1993) 187-216 \,\,\,\, MR1239560

\bibitem{segal} {\bf G B Segal,}  \emph{Classifying spaces and spectral sequences,} Publ. Math. Inst. des
Hautes Etudes Scient. (Paris) 34 (1968) 105-112  \,\,\,\, MR0232393

\bibitem{segal74} {\bf G B Segal,}  \emph{Categories and cohomology theories,} Topology  13 (1974)
293-312 \,\,\,\, MR0353298

\bibitem{sim} {\bf C Simson,}  \emph{A closed model structure for $n$-categories, internal $hom$,
 $n$-staks and generalized Seifert-Van Kampen} (arXiv:9704006)

\bibitem{street72} {\bf R Street,}  \emph{Two constructions on lax functors,}
Cahiers Top. G\'eom. Diff. 13 (1972) 217-264  \,\,\,\, MR0347936

\bibitem{street2} {\bf  R Street,}  \emph{ The algebra of oriented simplexes,}  J.  Pure Appl. Algebra
 (3) 49  (1987) 283-335  \,\,\,\, MR0920944

\bibitem{street} {\bf R Street,}  \emph{Categorical structures,} in Handbook of Algebra, Vol. 1,
 North-Holland, Amsterdam, (1996) 529-577 \,\,\,\, MR1421811


\bibitem{tam} {\bf Z Tamsamani,}   \emph{ Sur des notions de $n$-cat\'{e}gorie et $n$-groupo\"{i}de non strictes
 via des ensembles multi-simpliciaux,} K-theory 16 (1999)  51-59 \,\,\,\, MR1673923

\bibitem{thomason} {\bf R W Thomason,}  \emph{Homotopy colimits in the category of small categories,} Math. Proc. Camb. Phil. Soc. (1) 85
(1979)91-109  \,\,\,\, MR0510404

\bibitem{til}  {\bf U Tillmann,} \emph{On the homotopy of the stable mapping class group,} Invent.
Math.  (2) 130  (1997) 257-275  \,\,\,\, MR1474157

\bibitem{til2}  {\bf U Tillmann,} \emph{ Discrete models for the category of Riemann surfaces,} Math. Proc. Camb.
Phil. Soc. (1) 121  (1997) 39-49  \,\,\,\, MR1418359

\bibitem{tonks} {\bf K Worytkiewicz, K Hess, P E  Parent, A Tonks,}  \emph{ A model structure \`{a} la Thomason on 2-Cat,}
 J.  Pure Appl. Algebra  (1) 208 (2007) 205-236  \,\,\,\, MR2269840

\end{thebibliography}
\end{document}